\documentclass[12pt]{article}  
\usepackage[francais]{babel}
\usepackage{amssymb,amsmath}
\usepackage[ansinew]{inputenc}
\usepackage[dvips]{graphicx}

\newtheorem{theo}{Th{\'e}or{\`e}me}
\newtheorem{lemme}{Lemme}
\newtheorem{prop}[lemme]{Proposition}
\newtheorem{cor}[lemme]{Corollaire}
\newtheorem{defi}[lemme]{D{\'e}finition}
\newtheorem{rem}[lemme]{Remarque}
\newtheorem{exemple}[lemme]{Exemple}
\def\R{{\bf R}}
\def\C{{\bf C}}
\def\Q{{\bf Q}}
\def\Z{{\bf Z}}
\def\N{{\bf N}}
\def\H{{\bf H}}
\def\O{{\bf O}}
\def\hr{\mathcal{H}_\rho}
\def\me{m(\emptyset)}
\def\qed{~\vrule height 5pt width 5pt depth -.1pt}
\def\preuve{\par\medskip{\bf Preuve.}\hskip5pt}
\def\floor#1{\lfloor #1\rfloor}
\def\nn{|\hskip-1.5pt|}

\title{Superrigidité géométrique et applications harmoniques}
\author{Pierre Pansu$^{1,2}$\footnote{$^{1}$ Univ Paris-Sud, Laboratoire de Mathématiques d'Orsay, Orsay, F-91405 ;\hfill\eject\indent\hskip7.8pt $^{2}$ CNRS, Orsay, F-91405.}
}
\begin{document}
\maketitle

\section{Introduction}

\subsection{De quoi s'agit-il ?}

Le terme \emph{superrigidité}, inauguré par G.D. Mostow lorsqu'il a lu la contribution de G.A. Margulis au congrès de Vancouver, \cite{Margulis-ICM}, désigne un phé\-no\-mène mis en évidence par ce dernier en 1974 : les re\-pré\-sen\-ta\-tions linéaires de dimension finie non unitaires des réseaux de certains groupes de Lie proviennent, par restriction, du groupe ambiant. La situation modèle est celle du réseau $Sl(n,\Z)$ du groupe de Lie $Sl(n,\R)$.

Par \emph{superrigidité géométrique}, on entend un ensemble de techniques permettant d'étendre les résultats de Margulis à des classes plus larges de groupes discrets, et éventuellement de passer des re\-pré\-sen\-ta\-tions linéaires de dimension finie à des actions sur des espaces plus généraux.

\subsection{Le théorème de superrigidité de Margulis (1974)}
\label{margulis}

\begin{theo}
\label{marg}
{\em (Margulis, \cite{Margulis-ICM}).}
Soient $G$, $H$ des groupes algébriques semi-simples sur des corps locaux, sans facteurs compacts. On suppose que $G$ a un rang réel $\geq 2$. Soit $\Gamma$ un réseau irréductible de $G$. 

Tout homomorphisme $\Gamma\to H$ dont l'image est non bornée et Zariski dense s'étend en un homomorphisme $G\to H$.
\end{theo}

Nous ne définirons pas tous les termes, renvoyant à la littérature classique, \cite{Borel-livre}, \cite{Margulis-livre}. La géométrie commence lorsqu'on voit les groupes algébriques $G$ et $H$ comme groupes d'isométries d'espaces métriques. Depuis E. Cartan \cite{Cartan}, on sait qu'il existe un dictionnaire entre les groupes de Lie semi-simples sur $\R$ ou $\C$ et les \emph{espaces symétriques} sans facteurs euclidiens. Une variété riemannienne est symétrique si pour chaque point $x$, la symétrie géodésique, qui renverse toutes les géodésiques passant par $x$, est une isométrie. Ce dictionnaire a été étendu par F. Bruhat et J. Tits \cite{Bruhat-Tits} au cas des corps locaux non archimédiens. La classe des espaces symétriques est remplacée par celle des \emph{immeubles euclidiens}. Ce dictionnaire permet de reformuler le résultat de Margulis.

\begin{theo}
\label{marggeom}
Soient $X$, $Y$ des espaces symétriques ou des immeubles de dimension finie, sans facteurs compacts. On suppose que $X$ est de rang $\geq 2$. Soit $\Gamma$ un groupe discret irréductible d'isométries de $X$ tel que $Vol(\Gamma\setminus X)<+\infty$.

Toute action isométrique réductive de $\Gamma$ sur $Y$ laisse stable ou bien un point, ou bien un sous-ensemble convexe de $Y$ qui est pluriisométrique à un produit de facteurs irréductibles de $X$, et sur lequel l'action se prolonge en une action isométrique d'un quotient de $Isom(X)$.
\end{theo}

Certains termes sont plus aisés à définir dans ce langage. Lorsque $X$ s'écrit non trivialement comme un produit riemannien, un groupe d'isométries de $X$ est dit \emph{irréductible} s'il ne contient aucun sous-groupe de type fini qui préserve la dé\-com\-po\-si\-tion en produit. Toujours lorsque $X$ s'écrit non trivialement comme un produit riemannien, on peut modifier la métrique en multipliant celle de chaque facteur par une constante différente. On obtient ainsi les espaces \emph{pluriisométriques} à $X$. Le \emph{rang} de $X$ est la dimension maximale d'un espace euclidien qu'on peut plonger isométriquement dans $X$. Une action est \emph{réductive} (M. Gromov dit \emph{stable}, N. Monod dit \emph{non évanescente}), si pour tous les éléments $g$ d'un système générateur fini, la fonction déplacement $y\mapsto d(y,gy)$ sur $Y$ tend vers l'infini lorsque $d(y,y_0 )$ tend vers l'infini.

\bigskip

La restriction sur le rang a pu être partiellement levée.

\begin{theo}
\label{corlette}
\emph{(Corlette \cite{Corlette}, Gromov-Schoen \cite{Gromov-Schoen}).} Les résultats précédents s'é\-ten\-dent au cas où $X$ est un espace hyperbolique quaternionien de dimension $>4$ ou le plan hyperbolique des octonions.
\end{theo}

\begin{rem}
\label{rangun}
En revanche, ils ne s'étendent pas aux autres espaces sy\-mé\-triques de rang un (les espaces hyperboliques réels $\R H^n$ et complexes $\C H^n$) ni aux immeubles de rang un.
\end{rem}
C'est particulièrement frappant pour le plan hyperbolique réel et les immeubles de rang un, qui possèdent des réseaux qui sont des groupes libres. Il existe en toutes dimensions des variétés hyperboliques réelles qui possèdent une involution isométrique dont le lieu des points fixes, une hypersurface totalement géodésique, sépare la variété en deux, \cite{Millson}. Les réseaux de $SO(n,1)$ correspondants admettent une décomposition non triviale en produit amalgamé $\Gamma=A*_{C}B$. Le centralisateur de $C\subset SO(n,1)\subset SO(n+1,1)$ contient un $SO(2)$ qui ne centralise pas $A$ ni $B$. Pour $t\in SO(2)$ assez petit, l'homomorphisme $\rho_t :\Gamma\to SO(n+1,1)$ qui est l'identité sur $A$ et la conjugaison par $t$ sur $B$ est Zariski dense, mais ne correspond à aucun homomorphisme $SO(n,1)\to SO(n+1,1)$.

\subsection{Arithméticité}
Margulis a montré qu'on pouvait déduire d'un théorème de rigidité le fait que les réseaux sont arithmétiques, \cite{Margulis-arith}.

\begin{theo}
\label{arithm}
Les réseaux irréductibles des groupes de Lie semi-simples $G$ autres que $PO(n,1)=Isom(\R H^n )$ et $PU(n,1)=Isom(\C H^n )$ sont \emph{arith\-mé\-tiques}, i.e. obtenus (au relèvement à un produit $G\times L$ où $L$ est compact et à commensurabilité près) comme le groupe des matrices entières dans une re\-pré\-sen\-ta\-tion linéaire de $G$ définie sur $\Q$.
\end{theo}

Comme les groupes algébriques sur $\Q$ ont été classifiés (J. Tits, \cite{Tits-classification}), on aboutit à une classification des réseaux à commensurabilité près. C'est sans doute la conséquence la plus frappante du théorème de superrigidité.

\bigskip

A nouveau, le théorème \ref{arithm} ne s'étend pas aux groupes de rang 1 restants. Les mêmes variétés hyperboliques à symétrie évoquées au paragraphe \ref{margulis} servent à faire des réseaux non arithmétiques en toutes dimensions, \cite{Gromov-Piatetski-Shapiro} : M. Gromov et I. Piatetski-Shapiro choisissent soigneusement deux telles variétés, telles que les hypersurfaces séparantes soient isométriques, et recollent les moitiés respectives. Des réseaux non arithmétiques dans $SU(2,1)$ et $SU(3,1)$ ont ete construits par G.D. Mostow, \cite{Mostow-nonarithm}.

\subsection{Espaces $CAT(0)$}

Dès 1925, E. Cartan a observé que les espaces symétriques de type non compact satisfont des inégalités métriques là où l'espace euclidien satisfait des identités. Il s'en est servi pour montrer que tout groupe compact d'isométries possède un point fixe, et par conséquent, que les sous-groupes compacts maximaux du groupe des isométries sont deux à deux conjugués.

\begin{defi}
\label{cat}
Soit $Y$ un espace métrique géodésique. Etant donné un triangle de côtés $a$, $b$ et $c$ dans $Y$, on construit le triangle de même côtés $a'=a$, $b'=b$ et $c'=c$ dans le plan euclidien. A un point $u$ du côté $b$ correspond un point $u'$ qui divise le côté $b'$ dans les mêmes proportions. On note $d$ (resp. $d'$) la distance de $u$ (resp. $u'$) au sommet opposé. On dit que $Y$ est $CAT(0)$ si pour tout triangle, $d'\geq d$.
\end{defi}

\begin{center}
\includegraphics[width=3in]{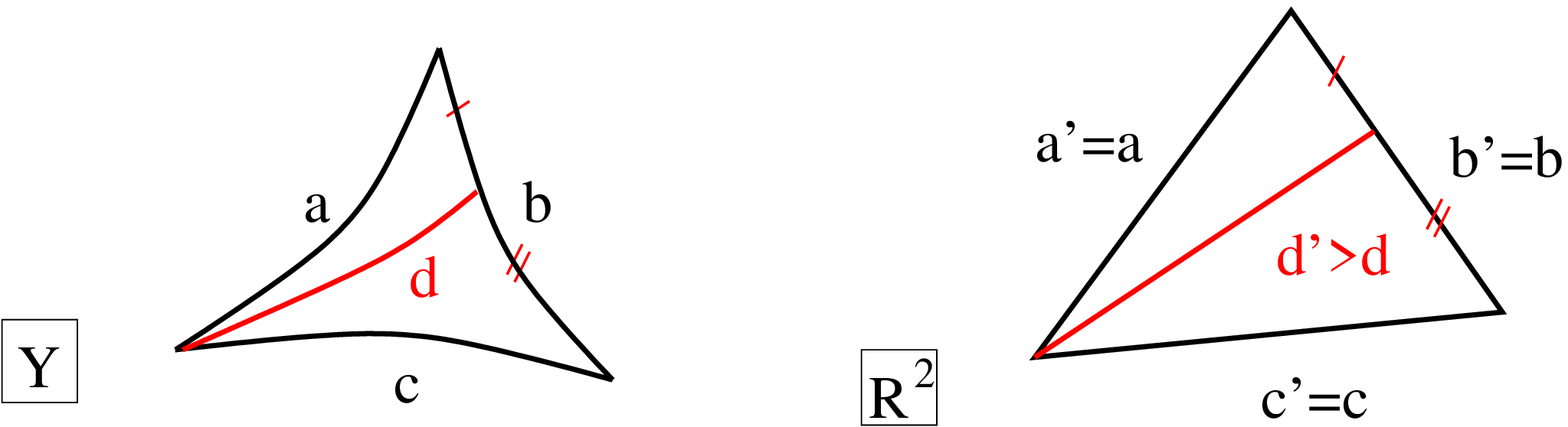}
\end{center}

L'observation de Cartan entraîne que les espaces symétriques de type non compact sont $CAT(0)$. De même, les arbres, les immeubles euclidiens sont géodésiques et $CAT(0)$, voir \cite{Bridson-Haefliger}.

Pour une variété riemannienne, être $CAT(0)$ est équivalent à être simplement connexe et à courbure sectionnelle négative ou nulle. Un produit d'espaces $CAT(0)$ est $CAT(0)$.

\subsection{Généralisation}
\label{rattaggi}

On peut baptiser superrigidité géométrique le programme suivant.

\medskip

\textbf{Question}. Soit $X$ un espace métrique. Trouver des conditions sur $X$ pour que pour tout groupe discret $\Gamma$ d'isométries de $X$, de covolume fini, toute action isométrique de $\Gamma$ sur un espace $CAT(0)$ géodésique et complet $Y$ ou bien possède un point fixe, ou bien laisse stable un sous-ensemble convexe de $Y$ pluriisométrique à $X$.

\begin{exemple}
Un théorème classique de V. Schroeder, \cite{Schroeder}, faisant suite à \cite{Gromoll-Wolf}, \cite{Lawson-Yau}, affirme que les groupes abéliens libres possèdent cette propriété de superrigidité géométrique, au moins pour les actions sur les variétés riemanniennes.
\end{exemple}
\begin{exemple}
Un résultat de N. Monod, \cite{Monod}, étendu par T. Gelander, A. Karlsson et G.A. Margulis, \cite{Gelander-Karlsson-Margulis}, prouve la superrigidité, sous une hypothèse de non évanescence, pour les réseaux uniformes irréductibles des produits de groupes localement compacts. 
\end{exemple}

Voici le type d'applications envisagées.

\begin{itemize}
  \item Expliquer en quoi les réseaux de $PU(n,1)$ ne sont pas superrigides.
  \item Etudier les re\-pré\-sen\-ta\-tions de dimension infinie des réseaux. 
  \item Etudier les actions sur des variétés compactes, au moyen de l'action induite sur un espace auxiliaire, comme l'espace des métriques riemanniennes mesurables sur la variété.
  \item Prouver que certains groupes ne sont pas \emph{linéaires}, i.e. isomorphes à des groupes de matrices.
\end{itemize}

Pour ces applications, on renvoie au survol \cite{Pansu1}, bien qu'il date un peu, et au texte de J. Maubon dans ce volume. Pour compléter ces survols, voir les trois références récentes \cite{Klingler}, \cite{Monod}, \cite{Gelander-Karlsson-Margulis}. 

Dans ces notes, on s'intéresse spécialement à la dernière, plus récente. 

\subsection{Groupes de type de re\-pré\-sen\-ta\-tion fini}

Il existe des groupes de présentation finie qui ne possèdent aucune re\-pré\-sen\-ta\-tion linéaire de dimension finie non triviale (par exemple, les groupes infinis simples, mais il y en a d'autres, voir \cite{Bass}). A mi-chemin entre ce phénomène et la non-linéarité (absence de re\-pré\-sen\-ta\-tions linéaire de dimension finie fidèles), Hyman Bass a introduit la notion suivante.

\begin{defi}
\emph{(H. Bass, \cite{Bass})}.
On dit qu'un groupe $\Gamma$ a un \emph{type de re\-pré\-sen\-ta\-tion fini} si toute re\-pré\-sen\-ta\-tion linéaire de dimension finie de $\Gamma$ a une image finie, i.e. factorise par un quotient fini de $\Gamma$.
\end{defi}

Le but de ces notes est d'expliquer comment une méthode, celle des applications harmoniques, inaugurée par J. Eells et J. Simons en 1964, mise en oeuvre dans le contexte de la superrigidité par Y.T. Siu en 1980, peut permettre de prouver qu'un groupe est de type de re\-pré\-sen\-ta\-tion fini.

\medskip

Je suis reconnaissant à M. Gromov, F. Haglund, M. Pichot et L. Silberman pour l'aide qu'ils m'ont apportée lors de la rédaction de ce texte. Celui-ci doit aussi beaucoup aux commentaires du referee.

\section{Superrigidité et finitude des re\-pré\-sen\-ta\-tions}

Voyons comment déduire la finitude des re\-pré\-sen\-ta\-tions d'une forme de superrigidité. Il s'agit d'un argument classique, dont on trouve les éléments dans \cite{Tits-libre}, \cite{Bass}, \cite{Margulis-arith}.

\begin{defi}
\label{deffsi}
On dit qu'un groupe $\Gamma$ possède la propriété FSI si toute action isomé\-tri\-que de $\Gamma$ sur un espace symétrique sans facteur compact ou un immeuble euclidien associé à $Sl(n)$ sur un corps local non archimédien possède un point fixe.
\end{defi}
Noter que l'espace euclidien est symétrique.

\begin{prop}
\label{FSI}
Soit $\Gamma$ un groupe de type fini. Si $\Gamma$ possède la propriété FSI, alors $\Gamma$ est de type de re\-pré\-sen\-ta\-tion fini.
\end{prop}

La preuve de la proposition \ref{FSI} occupe le reste de cette section. Voici un schéma de la preuve.

\bigskip

Soit $\Gamma$ un groupe qui possède la propriété FSI.
Etant donnée une re\-pré\-sen\-ta\-tion linéaire indécomposable de dimension finie, i.e. un homomorphisme $h:\Gamma\to Gl(n,\C)$, on montre successivement que
\begin{itemize}
  \item $h$ est localement rigide, i.e. un homomorphisme $h':\Gamma\to Gl(n,\C)$ suffisamment voisin de $h$ est conjugué à $h$ ;
  \item à conjugaison près, $h$ est à valeurs dans $Gl(n,\bar{\Q})$ ;
  \item après extension des scalaires et restriction à un sous-groupe d'indice fini de $\Gamma$, on obtient un homomorphisme $\tilde{h}$ à valeurs dans $Gl(nd,\Z)$ ;
  \item $\tilde{h}$ est à valeurs dans un sous-groupe compact de $Gl(nd,\C)$.
\end{itemize}
Il en résulte que $h(\Gamma)$ est fini.

\subsection{Cohomologie}

\begin{lemme}
\label{rigcohom}
Soit $\rho$ une re\-pré\-sen\-ta\-tion unitaire d'un groupe $\Gamma$ sur un espace de Hilbert $\hr$. Il y a une bijection entre l'ensemble des actions isométriques affines de $\Gamma$ sur $\hr$ (à conjugaison près par une translation), de partie linéaire $\rho$, et l'espace de cohomologie $H^1 (\Gamma,\rho)$. La classe associée à une action s'annule si et seulement si l'action possède un point fixe. 
\end{lemme}

\preuve
On suppose donnée une action isométrique affine de $\Gamma$ sur $\hr$, de partie linéaire $\rho$. Etant donné un point $v_0 \in\hr$, on pose, pour $g\in \Gamma$, $\eta_0 (g)=gv_0 -v_0$. Alors $\eta_0$ est un 1-cocycle sur $\Gamma$ à valeurs dans $\hr$. Sa classe de cohomologie ne dépend pas du choix de $v_0$. En fait, chaque cocycle $\eta_1$ cohomologue à $\eta_0$ correspond à un autre choix $v_1$ du point base. Par conséquent, si l'action possède un point fixe, $\eta$ est cohomologue à zéro et réciproquement.

Si une seconde action $(g,v)\mapsto g\cdot v$ donne le même cocycle, alors $gv_0 -v_0=\eta(g)=g\cdot v_0 -v_0$, donc $gv_0=g\cdot v_0$, et les deux actions, qui ont même partie linéaire, coïncident.

Inversement, soit $\eta$ un 1-cocycle sur $\Gamma$ à valeurs dans le $\Gamma$-module $\hr$. On construit une action (à gauche) isométrique affine de $\Gamma$ sur $\hr$ en posant, pour $v\in\hr$ et $g\in\Gamma$,
\begin{eqnarray*}
gv=\rho(g)(v)+\eta(g)
\end{eqnarray*}
Le cocycle $\eta_0$ qui lui est associé, lorsqu'on choisit comme point base $v_0 =0$, est $\eta$, par construction.\qed

\medskip

On conclut que, pour un groupe ayant la propriété $FSI$, $H^1 (\Gamma,\rho)=0$ pour toute re\-pré\-sen\-ta\-tion unitaire $\rho$ de dimension finie.

Lorsque $\rho$ est la re\-pré\-sen\-ta\-tion triviale, on trouve que $H^1 (\Gamma,\C)=0$. En particulier, tout homomorphisme de $\Gamma$ dans $\C$ est trivial. Aussi, l'abélianisé $\Gamma^{ab}=[\Gamma,\Gamma]$ est fini.

\subsection{Caractères}

Si $\Gamma$ est un groupe de type fini ayant la propriété $FSI$, alors tout homomorphisme $h:\Gamma\to Gl(1,\C)$ est d'image finie. En effet, un tel homomorphisme factorise par l'abélianisé $\Gamma^{ab}$. On conclut que pour tout homomorphisme $h:\Gamma\to Gl(n,\C)$, il existe un sous-groupe $\Gamma'$ d'indice fini dans $\Gamma$ tel que $h(\Gamma')\subset Sl(n,\C)$.

Dans la suite, on fera comme si $\Gamma'=\Gamma$.

\subsection{Unitarisabilité}

Soit $h:\Gamma\to Sl(n,\C)$ un homomorphisme. On fait agir $\Gamma$ via $h$ sur l'espace symétrique $Sl(n,\C)/SU(n)$. Par $FSI$, il existe un point fixe, c'est un sous-groupe compact conjugué de $SU(n)$ qui contient $h(\Gamma)$. Autrement dit, $h$ est conjuguée à une re\-pré\-sen\-ta\-tion unitaire. 

Le lemme \ref{rigcohom} s'applique donc à toutes les re\-pré\-sen\-ta\-tions linéaires de dimension finie de $\Gamma$.

\subsection{Rigidité locale}

On s'intéresse à l'espace $Rep=Hom(\Gamma,Gl(n,\C))/Gl(n,\C)$ des classes d'é\-qui\-va\-lence de re\-pré\-sen\-ta\-tions unitaires de $\Gamma$. Plus précisément, à sa structure locale au voisinage d'une re\-pré\-sen\-ta\-tion $h$. Notons $C(h)$ le centralisateur de $h(\Gamma)$ dans $Gl(n,\C)$. Depuis A. Weil \cite{W2}, on sait construire dans un produit de copies de $Gl(n,\C)$ une sous-variété analytique $C(h)$-invariante $PreRep$ dont l'espace tangent s'identifie à $H^1 (\Gamma,\rho)$ où $\rho=Ad\circ h$, et une application $C(h)$-é\-qui\-va\-rian\-te $\zeta:Prerep\to H^2 (\Gamma,\rho)$ telles que $Rep$ s'identifie au quotient par $C(h)$ du sous-ensemble analytique $C(h)$-invariant $\zeta^{-1}(0)$ de $PreRep$ (le couple $(Prerep,\zeta)$ s'appelle parfois modèle de Kuranishi, en référence à \cite{Kuranishi}). L'annulation de $H^1 (\Gamma,\rho)$ entraîne donc que la classe d'équivalence de $h$ est un point isolé de $Rep$.

\medskip

Démontrons directement un cas particulier de ce théorème. L'argument con\-tour\-ne l'obstacle du passage au quotient, qui donne parfois des espaces non séparés.

\begin{lemme}
\label{weil}
Soit $\Gamma$ un groupe de type fini. Soit $H$ un groupe de Lie. Soit $h:\Gamma\to H$ un homomorphisme dont l'image a un centralisateur discret dans $H$. On suppose que $H^1 (\Gamma,Ad\circ h)=0$. Alors tout homomorphisme $h':\Gamma\to H$ suffisamment proche de $h$ est conjugué à $h$.
\end{lemme}

\preuve
On considère $Hom(\Gamma,H)$ comme le sous-ensemble de $H^{\Gamma}$ défini par l'équation $\Phi(f)=1$ où $\Phi:H^{\Gamma}\to H^{\Gamma\times\Gamma}$, est définie comme suit. Etant donnée une fonction $f:\Gamma\to H$, $\Phi(f)$ est la fonction sur $\Gamma\times\Gamma$ définie par
\begin{eqnarray*}
\Phi(f)(g,g')=f(gg')^{-1}f(g)f(g').
\end{eqnarray*}
En fait, par noethérianité, un nombre fini d'équations suffisent.

Ramenons tous les vecteurs tangents à $H$ dans l'algèbre de Lie $\mathcal{H}$ par translation à gauche. Alors l'espace tangent en $h$ à $H^{\Gamma}$ s'identifie aux fonctions de $\Gamma$ dans $\mathcal{H}$, i.e. aux 1-cochaînes sur $\Gamma$ à valeurs dans $\mathcal{H}$, et l'espace tangent en 1 à $H^{\Gamma\times\Gamma}$ aux 2-cochaînes sur $\Gamma$ à valeurs dans $\mathcal{H}$. La différentielle de $\Phi$ en $h$ s'identifie (au signe près) au cobord relatif à la re\-pré\-sen\-ta\-tion $\rho=Ad\circ h$ : si $\eta\in C^1 (\Gamma,\hr)$, $h\eta\in T_h H^{\Gamma}$ et
\begin{eqnarray*}
D_h \Phi(h\eta)=-d\eta\in C^2 (\Gamma,\hr).
\end{eqnarray*}

L'action de $H$ par conjugaison sur $H^{\Gamma}$ définit une application $\Psi_h:H\to H^{\Gamma}$ comme suit. Etant donné $k\in H$, $\Psi_h (k)$ est la fonction sur $\Gamma$ définie par
\begin{eqnarray*}
\Psi_h (k)(g)=Ad_k (h(g))=k^{-1}h(g)k. 
\end{eqnarray*}
Comme $h$ est un homomorphisme, $\Phi\circ\Psi_h \equiv 1$.
Identifions $T_1 H=\mathcal{H}$ à l'espace des 0-cochaînes à valeurs dans $\mathcal{H}$. Alors la différentielle de $\Psi_h$ en 1 s'identifie (au signe près) au cobord : si $v\in \mathcal{H}$, qu'on voit comme une 0-cochaîne $c$,
\begin{eqnarray*}
D_1 \Psi_h (v)=-dc\in C^1 (\Gamma,\hr).
\end{eqnarray*}

On remarque que $\mathrm{ker}(d)=H^0 (\Gamma,\hr)=\{v\in\mathcal{H}\,|\forall g,\,\,Ad_{h(g)}(v)=v\}$ est l'algèbre de Lie du centralisateur de $h(\Gamma)$ dans $H$. Par hypothèse, cette algèbre de Lie est nulle. On a supposé aussi que
$H^1 (\Gamma,\hr)=0$. Par con\-sé\-quent, $D_1 \Psi_h$ est injective, et son image coïncide avec le noyau de $D_h \Phi$. Cette propriété est ouverte sur $Hom(\Gamma,H)$. Elle reste donc vraie pour $h'$ proche de $h$. Autrement dit, $\Phi$ est de rang constant sur $\Phi^{-1}(1)$ au voisinage de $h$. Par conséquent, $\Phi^{-1}(1)$ est une variété au voisinage de $h$, dont la dimension est égale à celle de l'orbite $\Psi_h (H)$, qu'elle contient. On conclut que, dans un voisinage de $h$, $Hom(\Gamma,H)$ coïncide avec l'orbite de $h$.

Dans le paragraphe précédent, on a fait comme si $H^{\Gamma}$ et $H^{\Gamma\times\Gamma}$ étaient des variétés de dimension finie. Voici pourquoi on peut le faire. Soit $S$ un système générateur fini symétrique de $\Gamma$. Choisissons, pour chaque élément $g\in\Gamma$, un mot $w(g)$ dans l'alphabet $S$ qui représente $g$. Toute fonction $f:S\to H$ (resp. $S\to\mathcal{H}$) se prolonge en une fonction $w(f):\Gamma\to H$ (resp. $\Gamma\to\mathcal{H}$), définie par $w(f)(g)=f(w(g))$. Si $h$ est un homomorphisme (resp. un 1-cocycle), il est uniquement déterminé par sa restriction à $S$ : $h=w(h_{|S})$. De plus, comme $Hom(\Gamma,H)$ est un ensemble analytique, pour vérifier qu'une fonction $f$ est un homomorphisme (resp. qu'une 1-cochaîne $\eta$ est un cocycle), il suffit de vérifier l'équation $f(gg')=f(g)f(g')$ (resp. $\eta(gg')=Ad_{h(g')}(\eta(g))+\eta(g')$)
pour un sous-ensemble fini $R\subset\Gamma\times\Gamma$ de couples $(g,g')$. 
On peut donc remplacer $\Phi:H^{\Gamma}\to H^{\Gamma\times\Gamma}$ par $\Phi_{S,R}=\pi_R \circ\Phi\circ w:H^S \to H^R$.\qed

\subsection{Algébricité}

Soit $\Gamma$ un groupe de type fini possédant la propriété $FSI$. Soit $h$ une re\-pré\-sen\-ta\-tion indécomposable de dimension finie de $\Gamma$. Par unitarisabilité, $h$ est irréductible. Par le lemme de Schur, le centralisateur de $h(\Gamma)$ dans le groupe spécial linéaire est réduit aux racines $n$-èmes de l'unité. D'après le lemme \ref{rigcohom}, $H^1 (\Gamma,Ad\circ h)=0$. Par rigidité locale, il existe un voisinage $V$ de $h$ dans $Hom(\Gamma,Sl(n,\C))$ tel que toute re\-pré\-sen\-ta\-tion contenue dans $V$ soit conjuguée à $h$. Comme $Hom(\Gamma,Sl(n,\C))$ est une sous-variété algébrique affine définie sur $\Q$ de $Gl(n,\C)^S$, $V$ contient un point défini sur $\bar{\Q}$. C'est un homomorphisme conjugué de $h$ dont l'image est contenue dans $Sl(n,\bar{\Q})$. On continue de le noter $h$.

\subsection{Intégralité}

Comme $\Gamma$ est de type fini, $h(\Gamma)$ est contenu dans $Sl(n,F)$ où $F$ est une extension de $\Q$ de degré fini $d$. Chaque élément de $F$ agit par multiplication sur $F$, d'où un plongement $F\to End_{\Q}(F)$, qui induit un plongement $Sl(n,F)\to End_{F}(F^n )\to End_{End_{\Q}(F)}(F^n )\simeq End_{\Q}(\Q^{nd})$, appelé \emph{extension des scalaires}.  On note $\tilde{h}$ sa composée avec $h$. L'image de $\tilde{h}$ est contenue dans $Sl(nd,\Q)$.

Soit $p$ un nombre premier. On fait agir $\Gamma$ via $\tilde{h}$ sur l'immeuble de $Sl(nd,\Q_p )$. Par superrigidité, il existe un point fixe, donc $\tilde{h}(\Gamma)$ est contenu dans un groupe compact, commensurable à $Sl(nd,\Z_p )$. Autrement dit, il existe un sous-groupe $\Gamma_p \subset\Gamma$ d'indice fini tel que $\tilde{h}(\Gamma_p )\subset Sl(nd,\Z_p )$.

L'étape précédente, appliquée à tous les diviseurs premiers des dé\-no\-mi\-na\-teurs des coefficients des matrices $\tilde{h}(s)$, $s\in S$, donne un sous-groupe $\Gamma' \subset\Gamma$ d'indice fini tel que $\tilde{h}(\Gamma' )\subset Sl(nd,\Z)$.

On fait agir $\Gamma$ via $\tilde{h}$ sur l'espace symétrique $Sl(nd,\C)/SU(nd)$. Par $FSI$, il existe un point fixe. Autrement dit, $\tilde{h}(\Gamma)$ est contenu dans un sous-groupe compact $K$, conjugué de $SU(nd)$. 

On conclut que $\tilde{h}(\Gamma' )\subset K\cap Sl(nd,\Z)$ est discret et relativement compact donc fini. Il en résulte que $h(\Gamma)$ est fini. Le cas général se ramène immédiatement au cas indécomposable.

\subsection{Groupes hyperboliques de type de re\-pré\-sen\-ta\-tion fini}

La superrigidité peut permettre de montrer que certains groupes ont un type de re\-pré\-sen\-ta\-tion fini. L'exemple suivant est dû à M. Kapovich, voir l'appendice de \cite{Kapovich}. Soit $\Gamma\subset G$ un réseau d'un des groupes de Lie simples de rang 1 qui sont superrigides, i.e. $G=Sp(n,1)$, $n\geq 2$, ou $F_{4}^{-20}$. Comme $\Gamma$ est hyperbolique au sens de Gromov, il admet des quotients infinis $\Gamma'=\Gamma/N$. Montrons que $\Gamma'$ possède la propriété FSI. Toute action isométrique de $\Gamma'$ sur un espace symétrique ou un immeuble de dimension finie $Y$ induit une action de $\Gamma$. D'après le Théorème \ref{corlette}, si cette action n'a pas de point fixe, elle laisse stable un sous-ensemble convexe $C\subset Y$ isométrique à un espace symétrique de rang un $X$, et l'action sur $C$ se prolonge à un quotient de $G$. Comme $G$ est simple, $G$ agit librement sur $C$, ce qui contredit le fait que le sous-groupe $N$ agit trivialement. On conclut que $\Gamma'$ fixe un point dans $Y$. Cela prouve que $\Gamma'$ a la propriété FSI, et donc un type de re\-pré\-sen\-ta\-tion fini.

En développant cet exemple, T. Kondo a montré que, dans l'espace des groupes marqués, il y a un $G_{\delta}$-dense de groupes de type fini qui sont hyperboliques et ont un type de re\-pré\-sen\-ta\-tion fini, voir \cite{Kondo}.

\section{Applications harmoniques é\-qui\-va\-rian\-tes}

Soit $\Gamma$ un groupe qui agit isométriquement sur des espaces $X$ et $Y$. Un point fixe de l'action, c'est la même chose qu'une application constante et é\-qui\-va\-rian\-te de $X$ dans $Y$. On va donc chercher une application constante parmi les applications é\-qui\-va\-rian\-tes $f$ de $X$ dans $Y$. Pour cela, on cherche à minimiser la ``variation quadratique'' de $f$, communément appelée {\em énergie}. Les applications qui minimisent l'énergie sont appelées {\em harmoniques}.

Cette stratégie remonte explicitement à J. Eells et J. Sampson en 1964. On peut même remonter plus loin. On a vu en \ref{rigcohom} qu'une action isométrique affine de $\Gamma$ sur un espace de Hilbert revient à se donner une re\-pré\-sen\-ta\-tion unitaire $\rho$ ainsi qu'une classe de cohomologie dans le $H^1 (\Gamma,\rho)$. Dans le cas où l'espace de Hilbert est de dimension 1 et la re\-pré\-sen\-ta\-tion triviale, une application é\-qui\-va\-rian\-te, à une constante près, c'est la même chose qu'une 1-forme différentielle fermée représentant la classe de cohomologie. L'idée de minimiser la norme $L^{2}$ dans une classe de cohomologie est à la base de la théorie de Hodge, née dans les années 30. Ce point de vue est développé dans \cite{Pansu1}. 

\subsection{Applications entre variétés}

Soit $Z$ une variété riemannienne compacte, $X$ son revêtement universel, $Y$ une variété riemannienne contractile, $f:X\to Y$ une application é\-qui\-va\-rian\-te. Alors la différentielle $df$ s'interprète comme une 1-forme dif\-fé\-ren\-tiel\-le sur $Z$ à valeurs dans un fibré vectoriel $E$ obtenu par passage au quotient du fibré induit $f^* TY$. La métrique et la connexion de Levi-Civita de $Y$ induisent une métrique et une connexion métrique $D$ sur le fibré $E$, et donc un opérateur $d^D$ sur les formes différentielles à valeurs dans $E$. En ce sens, $df$ est fermée,
\begin{eqnarray*}
d^D df =0.
\end{eqnarray*}

On définit l'\emph{énergie} de $f$ comme la moitié de la norme $L^2$ de $df$,
\begin{eqnarray*}
E(f)=\frac{1}{2}\int_{M}|df|^{2}.
\end{eqnarray*}
Les points critiques de l'énergie, appelés \emph{applications harmoniques}, satisfont l'é\-qua\-tion d'Euler-\-Lagrange
\begin{eqnarray*}
(d^{D})^{*}df=0,
\end{eqnarray*}
où $(d^{D})^{*}$ est l'adjoint du $d^D =D$ sur les 0-formes à valeurs dans $E$. On a $(d^{D})^{*}=tr\circ D$. Remarquer que les applications totalement géodésiques $X\to Y$ sont caractérisées par $Ddf=0$, et sont donc harmoniques. La superrigidité relativement à une classe de variétés riemanniennes se ramène donc
\begin{enumerate}
  \item à un théorème d'existence d'applications harmoniques é\-qui\-va\-rian\-tes,
  \item et un théorème d'annulation qui garantit que
  \begin{eqnarray*}
(d^{D})^{*}df=0 \Rightarrow Ddf=0.
\end{eqnarray*}
\end{enumerate}

\subsection{Résultats utilisant des applications harmoniques dé\-fi\-nies sur des variétés}

Ce programme en deux temps a été inauguré par J. Eells et J. Sampson en 1964. Dans \cite{Eells-Sampson}, ils montrent l'existence d'applications harmoniques é\-qui\-va\-rian\-tes dans le cas où $Y$ est à courbure sectionnelle négative ou nulle et ou $\Gamma=\pi_1 (Z)$ agit librement sur $Y$ avec quotient compact. De plus, une version à valeurs vectorielles de la formule de Bochner donne 
\begin{eqnarray*}
(d^{D})^{*}df=0 \Rightarrow Ddf=0.
\end{eqnarray*}
pourvu que la courbure de Ricci de $Z$ soit positive ou nulle. 

Cela entraîne la superrigidité (relativement à la classe des variétés riemanniennes à courbure sectionnelle négative ou nulle) des groupes abéliens de type fini.

\bigskip

En 1980, Y.T. Siu a utilisé une version vectorielle d'une formule spécifique aux variétés kählériennes, pour analyser les applications harmoniques définies sur des variétés kählériennes, sans hypothèse de signe de courbure, à valeurs dans des variétés complexes. J. Sampson a observé qu'il n'était pas indispensable que l'espace d'arrivée soit complexe. A sa suite, K. Corlette a découvert en 1990 un analogue quaternion-kählérien de la formule de Siu-Sampson, et prouvé la superrigidité (relativement à la classe des variétés riemanniennes à opérateur de courbure négatif ou nul) des réseaux de $H^{n}_{\H}$ et $H^{2}_{\O}$.

\bigskip

En 1994, M. Gromov et R. Schoen ont pu appliquer la formule de Corlette aux applications harmoniques à valeurs dans un immeuble. Pour cela, ils démontrent un résultat de régularité très fort : hors d'un lieu singulier de codimension au moins 2, l'application factorise par un plongement isométrique d'un espace euclidien dans l'immeuble. Combiné avec celui de Corlette, ce résultat constitue le théorème \ref{corlette}.

\bigskip

En 1993, en englobant la formule de Corlette dans une famille de formules dues à Y. Matsushima, N. Mok, Y.T. Siu and S.K. Yeung ont prouvé la superrigidité des réseaux irréductibles uniformes de rang supérieur (relativement à la classe des variétés riemanniennes à courbure sectionnelle négative ou nulle).

\bigskip

Si l'existence d'applications harmoniques é\-qui\-va\-rian\-tes d'une variété vers un espace CAT(0) général, sous diverses hypothèses de réductivité (ou stabilité ou non évanescence) est un problème résolu, voir \cite{Jost}, \cite{Korevaar-Schoen}, la difficulté de leur appliquer un théorème d'annulation limite actuellement l'extension de la méthode.

\section{Applications harmoniques combinatoires}

Ici, on décrit, en suivant M.T. Wang, \cite{Wang1}, \cite{Wang2} et H. Izeki et S. Nayatani \cite{Izeki-Nayatani}, une version discrète, élémentaire, de la notion d'application harmonique. Dans le cas particulier des actions isométriques affines sur des espaces de Hilbert, cette version discrète a, elle aussi, un glorieux passé (H. Garland, 1972, dans un langage cohomologique).

\subsection{Energie et applications harmoniques}
\label{sfixe}

Désormais, on suppose que l'espace $X$ est le revêtement universel d'un complexe simplicial fini dont $\Gamma$ est le groupe fondamental. On va définir une énergie pour les applications é\-qui\-va\-rian\-tes de $X$ 

\begin{defi}
\label{defpoids}
Soit $C$ un complexe simplicial. Un \emph{poids} sur $C$ est une fonction positive $m$ sur les simplexes telle que, pour tout $k$-simplexe $\sigma$, $m(\sigma)$ est égale à la somme des poids des $k+1$-simplexes qui contiennent $\sigma$,
\begin{eqnarray*}
m(\sigma)=\sum_{\sigma\subset\tau} m(\tau).
\end{eqnarray*}
\end{defi}

\begin{exemple}
\label{expoids}
Si $C$ est de dimension 2, on pose $m=1$ sur les 2-simplexes, et on propage aux 1-simplexes puis aux 0-simplexes en respectant la propriété de poids, et on normalise par le poids total.
\end{exemple}

\begin{defi}
\label{defdist}
Soit $(C,m)$ un complexe simplicial pondéré fini. On note $C^0$ l'ensemble de ses sommets, $C^1$ l'ensemble de ses arêtes. Soit $Y$ un espace métrique. Sur l'espace des fonctions de $C$ dans $Y$, on met la distance ``produit pondéré'' 
\begin{eqnarray*}
d(g,g')=(\sum_{c\in C^0 }m(c)d(g(c),g'(c))^{2})^{1/2}.
\end{eqnarray*}
\end{defi}

\begin{defi}
\label{defenerg}
Etant donnés un complexe simplicial pondéré fini $(C,m)$ et une application $f:C^0 \to Y$ à valeurs dans un espace métrique, l'\emph{énergie} de $f$ est
\begin{eqnarray*}
E(f)&=&\sum_{e\in C^1}m(e)d(f(ori(e)),f(ext(e)))^2 \\
&=&\frac{1}{2}\sum_{c\in C^0}\sum_{c'\sim c}m(c,c')d(f(c),f(c'))^2 ,
\end{eqnarray*}
où $c\sim c'$ signifie que $c$ et $c'$ sont les extrémités d'une même arête.
\end{defi}

Voici une version é\-qui\-va\-rian\-te des définitions précédentes.
\begin{defi}
\label{defenergequiv}
Soit $Z$ un complexe simplicial fini, de groupe fondamental $\Gamma$, soit $X$ son revêtement universel. On se donne un poids $m$ sur $Z$. Soit $Y$ un espace métrique sur lequel $\Gamma$ agit isométriquement. Sur l'espace des applications é\-qui\-va\-rian\-tes $X\to Y$, on met la distance
\begin{eqnarray*}
d(f,f')=(\sum_{z\in Z^0}m(z)d(f(c),f'(c))^{2})^{1/2},
\end{eqnarray*}
où $c$ est un sommet de $X$ au-dessus de $z$. L'\emph{énergie} de $f$ est
\begin{eqnarray*}
E(f)&=&\sum_{e\in Z^1}m(e)d(f(ori(\tilde{e})),f(extr(\tilde{e}))^2 \\
&=&\frac{1}{2}\sum_{z\in Z^0}\sum_{z'\sim z} m(z,z')d(f(c),f(c'))^2 ,
\end{eqnarray*}
où $\tilde{e}$ (resp. $(c,c')$) est une arête de $X$ au-dessus de l'arête $e$ (resp. $(z,z'))$ de $Z$.
\end{defi}

\begin{rem}
\label{remequiv}
Si $Y$ est $CAT(0)$, l'espace $Y^{\Gamma\setminus X}$ des applications é\-qui\-va\-rian\-tes $X\to Y$ est lui-même $CAT(0)$, et l'énergie est une fonction convexe sur $Y^{\Gamma\setminus X}$.
\end{rem}
En effet, $t\mapsto f_t$ est une géodésique parcourue à vitesse constante de $Y^{\Gamma\setminus X}$ si et seulement si pour tout $x\in X$, $t\mapsto f_t (x)$ est une géodésique parcourue à vitesse constante de $Y$. Comme le caractère CAT(0) de $Y$ s'exprime par une inégalité entre carrés de distances, elle s'intègre en la propriété CAT(0) pour $Y^{\Gamma\setminus X}$. D'autre part, dans un espace $CAT(0)$, étant données deux géodésiques $t\mapsto s(t)$ et $t\mapsto s'(t)$ parcourues à vitesse constante, la fonction $(t,t')\mapsto d(s(t),s'(t'))$ est convexe sur $\R^2$.

\begin{defi}
\label{defharm}
Une application é\-qui\-va\-rian\-te $f:X\to Y$ est dite {\em harmonique} si elle minimise localement l'énergie $E$, i.e. si un petit déplacement équivariant de l'image d'une orbite de sommets ne diminue pas l'énergie.
\end{defi}

Autrement dit, on interprète $X$ comme un arrangement périodique de ressorts, et $f$ comme une déformation périodique des ressorts (voir figure). $E$ représente l'énergie potentielle élastique par maille. Une application $f$ est harmonique si elle représente une position d'équilibre.

\begin{center}
\includegraphics[width=3in]{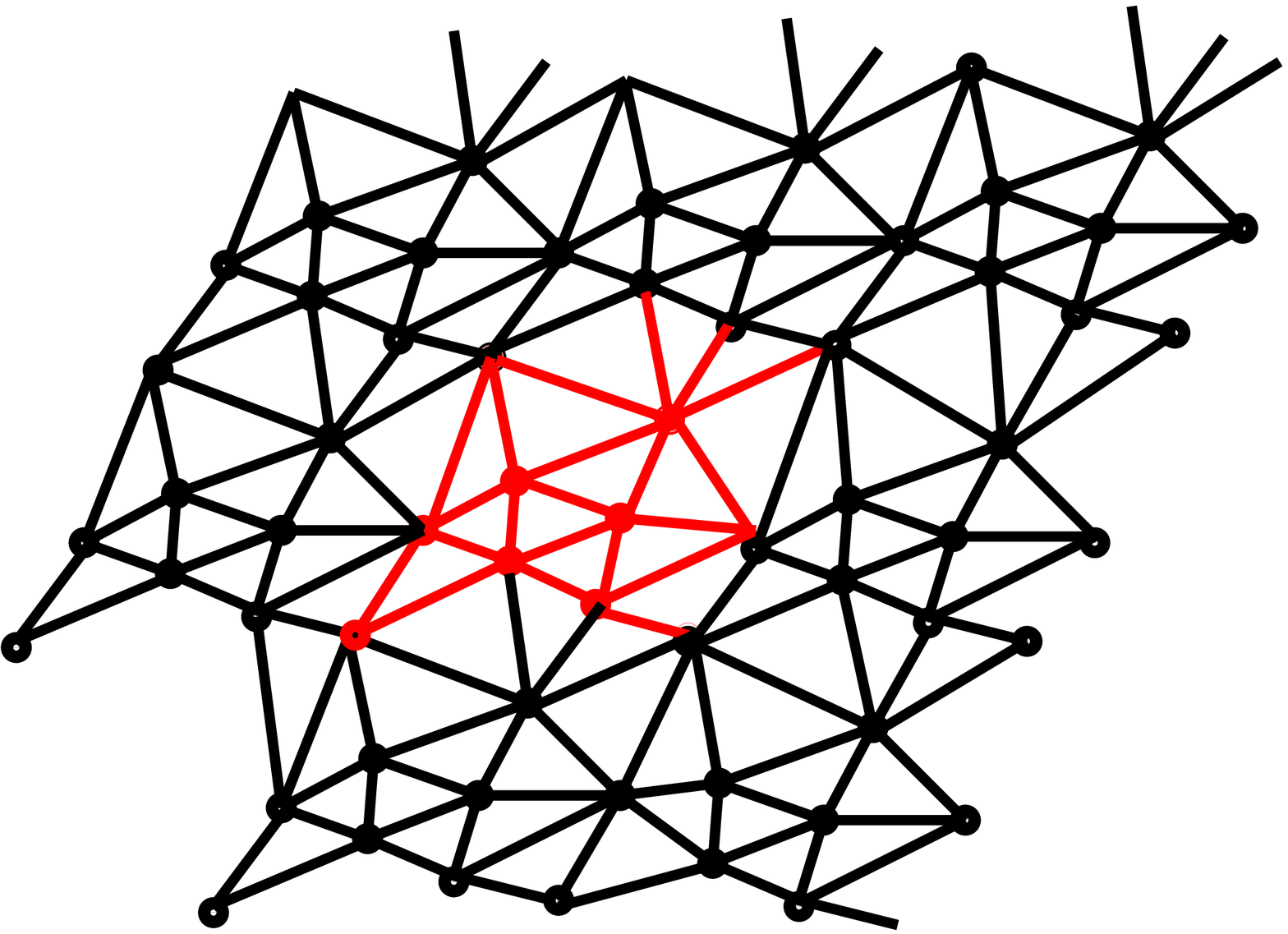}
\end{center}

\subsection{Barycentre, lien et équation des applications harmoni\-ques}

Dans l'espace euclidien, le carré $u$ de la distance à un point est une fonction uniformément strictement convexe : le long d'une géodésique $t\mapsto y(t)$ parcourue à vitesse 1, $u(y(t))-t^{2}$ est convexe. Par comparaison, cette propriété s'étend aux espaces géodésiques $CAT(0)$. Elle se généralise aux combinaisons à coefficients positifs de carrés de distances à des points, ce qui permet de montrer que ces fonctions atteignent leur minimum. 

\begin{lemme}
\label{lembar}
Soit $Y$ un espace métrique $CAT(0)$ complet, soit $\mu$ une mesure de probabilité sur $Y$. Posons, pour $y\in Y$, $u(y)=d(\mu,y)^2 =\int_Y d(z,y)^2 \,d\mu(z)$. Alors 
\begin{eqnarray*}
diam(\{u\leq \epsilon^2 +\inf u\})\leq 2\epsilon.
\end{eqnarray*}
Par conséquent, $u$ atteint son minimum sur $Y$, en un unique point $y_0$, et on a, pour tout $y\in Y$,
\begin{eqnarray*}
u(y)\geq u(y_0 )+d(y,y_0 )^2 .
\end{eqnarray*}
\end{lemme}

\preuve
Soit $t\mapsto s(t)$, $t\in[-L,L]$, un segment géodésique de vitesse 1, reliant deux points où $u\leq \epsilon^2 +\inf u$. La fonction $t\mapsto v(t)=u(t)-\inf u -\,t^2$ est convexe, donc $0\leq v(0)\leq \frac{1}{2}(v(-L)+v(L))\leq\epsilon^2 -\,L^2$, donc $L\leq\epsilon$.\qed

\medskip

Etant donné un complexe simplicial pondéré fini $(C,m)$ et une application $g:C^0 \to Y$, on applique ce lemme à la mesure image par $g$ de la mesure de probabilité $m$.

\begin{defi}
\label{defbar}
Soit $(C,m)$ un complexe simplicial pondéré fini, $Y$ un espace mé\-tri\-que $CAT(0)$ complet, $g:C^0 \to Y$ une application. L'unique point de $Y$ où la fonction $y\mapsto d(g,y)^{2}=d(g_* (m/\me),y)$ atteint son minimum est appelé {\em barycentre} de $g$.
\end{defi}

\begin{exemple}
\label{exbar}
Lorsque $Y$ est un espace de Hilbert, le barycentre métrique coïncide avec le barycentre affine $bar(g)=\frac{1}{\me}\sum_c m(c)g(c)$.
\end{exemple}

Etant donné un complexe simplicial pondéré $(X,m)$, l'{\em étoile} d'un sommet $x$ est la réunion des simplexes de $X$ qui contiennent $x$, et le {\em lien} de $x$ est la réunion des simplexes de l'étoile qui ne contiennent pas $x$.
 
\begin{center}
\includegraphics[width=2in]{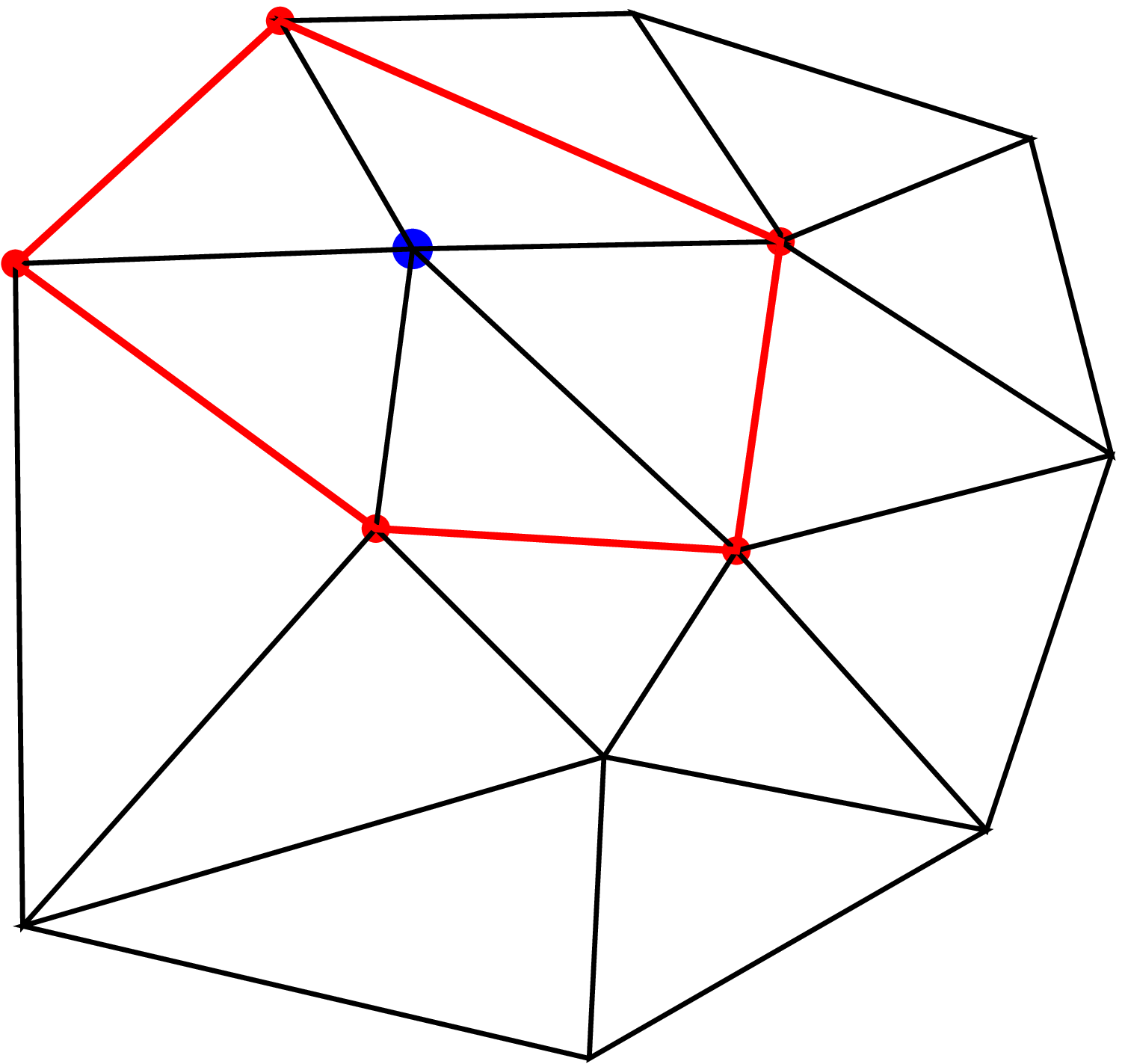}
\end{center}

Le lien de $x$ hérite d'une pondération : si $\sigma$ est un simplexe du lien, $x\cup\sigma$ est un simplexe de $X$, et on pose $m_{lien}(\sigma)=m(x\cup\sigma)$.

\begin{prop}
\label{eqharm}
Soit $(Z,m)$ un complexe simplicial pondéré fini de groupe fondamental $\Gamma$, soit $X$ son revêtement universel. Soit $Y$ un espace métrique $CAT(0)$. Une application é\-qui\-va\-rian\-te $f:X\to Y$ est harmonique si et seulement si pour tout sommet $x\in X$, $f(x)$ coïncide avec le barycentre de la restriction de $f$ au lien de $x$,
\begin{eqnarray*}
f(x)=bar(f_{|lien(x)}).
\end{eqnarray*}
\end{prop}

\preuve
En effet, si on déplace de façon é\-qui\-va\-rian\-te l'image de l'orbite d'un sommet $x$ de $X$, de sorte que $f(x)$ est remplacé par $y\in Y$, l'énergie de $f$ change de
\begin{eqnarray*}
\sum_{x'\in lien(x)}m(x,x')(d(f(x'),y)^{2}-d(f(x'),f(x))^{2}).
\end{eqnarray*} 
Ce terme est positif ou nul pour tout $y$ si et seulement si $f(x)=bar(f_{|lien(x)})$.\qed

\begin{exemple}
\label{exharm}
Lorsque $Y=\R$, on peut interpréter $X$ comme un circuit électrique, une fonction $f$ comme le potentiel électrique aux sommets du circuit, le nombre $m(x,x')$ comme la conductance de la branche reliant $x$ à $x'$, $m(x,x')(f(x')-f(x))$ comme le courant qui passe dans cette branche. La loi de Kirchhoff énonce alors que $f$ est harmonique.
\end{exemple}

\subsection{Flot du gradient}

D'après U. Mayer, \cite{Mayer}, on peut définir le gradient d'une fonction conve\-xe continue sur un espace $CAT(0)$ complet. Dans le cas de la fonction énergie sur l'espace des applications é\-qui\-va\-rian\-tes de $X$ vers un espace $CAT(0)$ complet $Y$ (elle est convexe, voir \ref{remequiv}), il est naturel de noter $\Delta f$ le gradient en $f$.

\begin{defi}
\label{deflap}
\emph{(H. Izeki et S. Nayatani, \cite{Izeki-Nayatani}).}
Soit $f:X\to Y$ une application é\-qui\-va\-rian\-te. Pour $x\in X$, on note $-\Delta f(x)\in T_{f(x)}Y$ le barycentre de la fonction $f_x =\pi\circ f_{|lien(x)}$ définie sur le lien de $x$, à valeurs dans le cône tangent à $Y$ en $f(x)$. On note $|-\Delta f(x)|$ la distance de $-\Delta f(x)$ au sommet du cône, et
\begin{eqnarray*}
|-\Delta f|^{2}=\sum_{z}m(z)|-\Delta f(x)|^{2}.
\end{eqnarray*}
\end{defi}
Remarquer que $f$ est harmonique si et seulement si $|-\Delta f|=0$.

\bigskip

U. Mayer a construit un flot de gradient pour toute fonction convexe continue sur un espace $CAT(0)$ complet. Cela s'applique à l'énergie sur l'espace des applications é\-qui\-va\-rian\-tes. On dispose donc, pour toute application initiale $f_0$, d'une ligne de gradient $t\mapsto f_t$, solution en un sens faible de l'équation différentielle
\begin{eqnarray*}
\frac{df_t}{dt}=-grad_{f_t}E=2(-\Delta f_t ).
\end{eqnarray*}
Alors
\begin{eqnarray*}
\frac{d}{dt}E(f_t )=-|grad_{f_t}E|^{2}=-4|-\Delta f_{t}|^{2}.
\end{eqnarray*}

\bigskip

Voici une autre conséquence de la convexité de l'énergie.

\begin{lemme}
\label{lemmayer}
Soient $f$, $g:X\to Y$ des applications é\-qui\-va\-rian\-tes. Alors
\begin{eqnarray*}
E(f)-E(g)\leq |-\Delta f|d(f,g).
\end{eqnarray*}
\end{lemme}

\preuve
Pour chaque $x\in X$, notons $t\mapsto g_t (x)$, $t\in[0,1]$, la géodésique de vitesse constante reliant $f(x)$ à $g(x)$. Alors $t\mapsto g_t$ est une géodésique parcourue à vitesse constante égale à $d(f,g)$ dans l'espace des applications é\-qui\-va\-rian\-tes. La fonction $t\mapsto E(g_t )$ est convexe, et sa dérivée en $t=0$ est comprise entre $-|-\Delta f(x)|d(f,g)$ et $|-\Delta f(x)|d(f,g)$, donc
\begin{eqnarray*}
E(g)=E(g_1 )\geq E(f)-|-\Delta f(x)|d(f,g).\qed
\end{eqnarray*}

\subsection{Ultralimites}

On ne peut pas toujours prouver la convergence du flot de Mayer vers une application constante. Néanmoins, une sous-suite converge vers une application harmonique, éventuellement à valeurs dans un autre espace, obtenu comme limite faible de dilatés de $Y$. Cette observation est due initialement à N. Mok, dans le cas où $Y$ est un espace de Hilbert, \cite{Mok} (ce résultat a été retrouvé indépendamment par Y. Shalom, \cite{Shalom}). 

Dans le cas où $Y$ n'est pas linéaire, la notion pertinente est celle d'ultralimite, classique en théorie des modèles, et popularisée par M. Gromov dans \cite{Gromov-AIIG}. 

Soit $Y_j$ une suite d'espaces métriques. On se donne une suite de points $y_j \in Y_j$ et un ultrafiltre non principal $\omega$ sur $\N$. On considère l'ensemble des suites $(x_j )_{j\in\N}$ de points $x_j \in Y_j$ telles que la distance $d(x_j ,y_j )$ reste bornée. On le munit de la semi-distance
$d((x_j ),(x'_j ))=\lim_{\omega}d(x_j ,x'_j )$. On quotiente par la relation $(x_j )\sim(x'_j )$ si $d((x_j ),(x'_j ))=0$. L'espace métrique obtenu est noté $\lim_{\omega}(Y_j )$. C'est une \emph{ultralimite} de la suite $Y_j$. Remarquer que lorsque les $Y_j$ sont des espaces géodésiques $CAT(0)$ complets, toute ultralimite des $Y_j$ est géodésique, $CAT(0)$ et complet, voir \cite{Gromov-AIIG}.

\begin{defi}
\label{cone}
On appelle \emph{cône asymptotique} de $(Y,d)$ toute ultralimite d'une suite $Y_j =(Y,d_j =r_j d,y_j )$ d'espaces homothétiques de $Y$, avec $r_j \to +\infty$ ou $y_j \to \infty$.
\end{defi}

\subsection{Existence d'applications harmoniques}
\label{existence}

\begin{theo}
\label{exharmlim}
\emph{(Enoncé par M. Gromov, \cite{Gromov-RWRG} paragraphe 3.6, voir aussi H. Iseki, T. Kondo and S. Nayatani, \cite{Izeki-Kondo-Nayatani}).} 
Soit $Z$ un complexe simplicial fini de groupe fondamental $\Gamma$. Soit $X$ son revêtement universel. Soit $Y$ un espace métrique gé\-o\-dé\-si\-que, $CAT(0)$ et complet. Soit $\rho$ une action isométrique de $\Gamma$ sur $Y$. Alors
\begin{enumerate}
  \item ou bien $\rho$ fixe un point ;
  \item ou bien il existe une application harmonique $\rho$-equivariante \emph{non cons\-tan\-te} de $X$ dans $Y$ ;
  \item ou bien il existe un cône asymptotique $Y_{\omega}$ de $Y$, une action isométrique $\rho_{\omega}$ de $\Gamma$ sur $Y_{\omega}$ et une application harmonique, $\rho_{\omega}$-é\-qui\-va\-rian\-te et \emph{non constante} de $X$ dans $Y_{\omega}$.
\end{enumerate}
De plus, s'il existe une constante $C$ telle que pour toute application é\-qui\-va\-rian\-te $f$,
\begin{eqnarray*}
E(f)\leq C|-\Delta f|^{2},
\end{eqnarray*}
alors on est dans le premier cas, $\rho$ fixe un point.
\end{theo}

\preuve
Soit $f:X\to Y$ une application é\-qui\-va\-rian\-te, soit $f_t$ le flot de Mayer. 

S'il existe une constante $C$ telle que pour tout $t$ assez grand,
\begin{eqnarray*}
E(f_t )\leq C|-\Delta f_t |^{2},
\end{eqnarray*}
on montre que $f_t$ converge vers une application constante. D'abord,
\begin{eqnarray*}
\frac{d}{dt}E(f_t )\leq -4C^{-1} E(f_t ),
\end{eqnarray*}
donc l'énergie décroît exponentiellement. 

Comme 
\begin{eqnarray*}
|\frac{df_t}{dt}|^2=-|grad_{f_t}E|^{2}=-\frac{d}{dt}E(f_t ),
\end{eqnarray*}
l'intégrale
\begin{eqnarray*}
\int_{0}^{+\infty}|\frac{df_t}{dt}|^2 E(f_{t})^{-1/2}\,dt=2 E(f_{0})^{1/2}
\end{eqnarray*}
est finie, d'où
\begin{eqnarray*}
\int_{0}^{+\infty}|\frac{df_t}{dt}|\,dt&\leq&
(\int_{0}^{+\infty}|\frac{df_t}{dt}|^2 E(f_{t})^{-1/2}\,dt)^{1/2}(\int_{0}^{+\infty}E(f_{t})^{1/2}\,dt)^{1/2}<\infty.
\end{eqnarray*}
Autrement dit, la courbe $t\mapsto f_t$ est de longueur finie. On conclut que $f_t$ converge vers une application constante, et $\rho$ fixe un point.

Sinon, il existe une suite $t_j$ tendant vers $+\infty$ telle que $\displaystyle \frac{|-\Delta f_{t_j} |^{2}}{E(f_{t_j} )}$ tend vers 0. Choisissons un ultrafiltre non principal $\omega$. Fixons une origine $o\in X$. Posons $y_j =f_{t_{j}}(o)$, $\displaystyle r_j =\frac{1}{\sqrt{E(f_{t_j} )}}$, $Y_j =(Y,d_j =r_j d,y_j )$. Vue comme une application $X\to Y_j$, l'énergie de $f_{t_{j}}$ vaut 1. Par conséquent, pour tous $x$, $x'\in X$, $d_j (f_{t_{j}}(x),f_{t_{j}}(x'))$ est borné. La suite de points $f_{t_{j}}(x)$ représente donc un point, noté $f_{\omega}(x)$, de $Y_{\omega}$. Si $\gamma\in\Gamma$, et si $w_{\omega}\in Y_{\omega}$ est représenté par une suite $(w_j )$ située à distance bornée de $y_j$, alors 
\begin{eqnarray*}
d_j (\gamma w_j ,y_j )=d_j (w_j ,f_{t_{j}}(\gamma^{-1}o))\leq d_j (w_j ,y_j )+d_j (f_{t_{j}}(o),f_{t_{j}}(\gamma^{-1}o))
\end{eqnarray*}
est borné, donc les points $\gamma w_j$ représentent un point de $Y_{\omega}$ noté $\gamma w_{\omega}$. On définit ainsi l'action $\rho_{\omega}$ de $\Gamma$ sur $Y_{\omega}$, et l'application $f_{\omega}$ est $\rho_{\omega}$-é\-qui\-va\-rian\-te.

Montrons que $f_{\omega}$ minimise l'énergie. Soit $g_{\omega}:X\to Y_{\omega}$ une autre application $\rho_{\omega}$-é\-qui\-va\-rian\-te. Pour chaque $z\in Z$, choisissons un représentant $x\in X$ et représentons $g_{\omega}(x)$ par une suite $g_{j}(x)$ située à distance bornée de $f_{t_{j}}(x)$. Prolongeons chaque $g_j$ en une application é\-qui\-va\-rian\-te $X\to Y_j$. Les distances $\ell^2$ $d_{j}(f_{t_{j}},g_{j})$ sont bornées. Le lemme \ref{lemmayer} donne
\begin{eqnarray*}
E_j (f_{t_{j}})-E_j (g_j )\leq |-\Delta_{j} f_{t_{j}}|d_{j}(f_{t_{j}},g_{j}).
\end{eqnarray*}
Or le changement de distance sur $Y$ a pour effet que
\begin{eqnarray*}
|-\Delta_{j} f_{t_{j}}|^{2}=\frac{|-\Delta f_{t_{j}}|^{2}}{E(f_{t_{j}})}
\end{eqnarray*}
tend vers 0, donc $E_{\omega}(f_{\omega})-E_{\omega}(g_{\omega})= \lim E_j (f_{t_{j}})-E_j (g_j )\leq 0$. Ceci prouve que $f_{\omega}$ est harmonique. Comme $E_{\omega}(f_{\omega})=1$, $f_{\omega}$ n'est pas constante. 

Si les suites $r_j$ et $y-j$ sont bornées, $Y_{\omega}$ est simplement homothétique de $Y$, et on a trouvé une application harmonique à valeurs dans $Y$. Sinon, $r_j$ tend vers $+\infty$, et $Y_{\omega}$ est un cône asymptotique de $Y$.\qed 

\section{Formule de Garland}

Il s'agit d'un mécanisme qui force les applications harmoniques à être constantes. Combiné avec le théorème \ref{exharmlim}, il entraîne immédiatement la propriété de point fixe.

Un invariant va jouer un rôle clé dans cette formule, c'est le \emph{bas du spectre} non linéaire introduit par M.T. Wang dans sa thèse.

\subsection{Bas du spectre}

\begin{defi}
\label{defbas}
Soit $(C,m)$ un complexe simplicial pondéré, $Y$ un espace métrique $CAT(0)$, $g:C^0 \to Y$ une application non constante. Le {\em quotient de Rayleigh} de $g$ est
\begin{eqnarray*}
RQ(g)=\frac{E(g)}{d(g,bar(g))^{2}}.
\end{eqnarray*}
(La distance $g$ à l'application constante $bar(g)$ est celle définie en \ref{defdist}).

Le {\em bas du spectre} de $(C,m)$ relativement à l'espace $Y$ est la borne in\-fé\-ri\-eu\-re des quotients de Rayleigh des applications non constantes de $C^0$ dans $Y$,
\begin{eqnarray*}
\lambda(C,m,Y)=\inf_{g:C^0 \to Y} RQ(g).
\end{eqnarray*}
\end{defi}

\begin{exemple}
\label{exbas}
Lorsque $Y=\R$, le bas du spectre coïncide avec la plus petite valeur propre non nulle du laplacien discret
\begin{eqnarray*}
g\mapsto \Delta g,\quad \textrm{où}\quad\Delta g(c)=\sum_{c'\sim c}m(c,c')(g(c)-g(c')).
\end{eqnarray*}
\end{exemple}
On donnera en section \ref{basspectre} des exemples de bas de spectres de graphes.

\begin{exemple}
\label{exbasprod}
$\lambda(C,m,Y_{1}\times Y_{2})=\min\{\lambda(C,m,Y_{1}),\lambda(C,m,Y_{2})\}$.
\end{exemple}
En effet, la distance sur le produit est la racine carrée de la somme des carrés des distances sur les facteurs. Une application $g:C^0 \to Y_{1}\times Y_{2}$ s'écrit $g=(g_{1},g_{2})$, le barycentre $bar(g)=(bar(g_{1}),bar(g_{2}))$, $RQ(g)=RQ(g_{1})$ (resp. $RQ(g_{2})$) si $g_{2}$ (resp. $g_{1}$) est constante, et est un barycentre de $RQ(g_{1})$ et de $RQ(g_{2})$ sinon.

De même, si $Y$ est un espace de Hilbert, alors pour tout complexe simplicial pondéré fini $(C,m)$, $\lambda(C,m,Y)=\lambda(C,m,\R)$.

\begin{lemme}
\label{baslim}
Soient $Y_j$ des espaces géodésiques $CAT(0)$ complets. Soit $Y$ une ultralimite des $Y_j$. Alors $Y$ est géodésique $CAT(0)$ complet, et pour tout graphe pondéré fini $(C,m)$,
\begin{eqnarray*}
\lambda(C,Y)\geq \limsup_{j\to \infty}\lambda(C,Y_j ).
\end{eqnarray*}
\end{lemme}

\preuve
Soit $\omega$ un ultrafiltre tel que $(Y,y)=\lim_{\omega}(Y_j ,y_j )$. Etant donnée $g:C\to Y$, choisissons pour chacun des points $g(c)\in Y$ une suite qui le représente. On la note $(g_j (c))_{j\in\N}$. Alors $\lim_{\omega}E(g_j )=E(g)$. La suite $(bar(g_j ))$ restant à distance bornée de l'image de $g_j$, elle représente un point $z\in Y$, et $\lim_{\omega}d(g_j ,bar(g_j ))^2=d(g,z)^2 \geq d(g,bar(g))^2$, donc \begin{eqnarray*}
\lim_{\omega}RQ(g_j )\leq RQ(g).
\end{eqnarray*}
On conclut que 
\begin{eqnarray*}
\limsup_{j\to \infty}\lambda(C,Y_j )\leq \lambda(C,Y).\qed
\end{eqnarray*}

\subsection{Formule de Garland}

Cette formule a été découverte initialement par H. Garland, \cite{Garland} : en fait, il s'agit d'une famille de formules qui s'appliquent aux cocycles harmoniques des quotients compacts des immeubles euclidiens, dans tous les degrés. A. Borel \cite{Borel-Garland} a su généraliser la première d'entre elles à des complexes simpliciaux quelconques. A. Zuk \cite{Zuk-96} a été le premier à en tirer partie en dimension infinie. La version non linéaire, due à Wang (\cite{Wang1} pour les variétés, \cite{Wang2} en général), a été retrouvée par M. Gromov.

\begin{theo}
\label{garland}
Soit $Z$ un complexe simplicial fini de groupe fondamental $\Gamma$. Soit $X$ son revêtement universel. Soit $Y$ un espace métrique géodésique, $CAT(0)$ et complet. Soit $f:X\to Y$ une application é\-qui\-va\-rian\-te. Pour $x\in X$, on note
\begin{eqnarray*}
ED(f,x)=\frac{1}{2}d(f_{|lien(x)},f(x))^2 ,
\end{eqnarray*}
de sorte que
\begin{eqnarray*}
E(f)=\sum_{z\in Z}ED(f,x).
\end{eqnarray*}
Alors, si $f$ est harmonique,
\begin{eqnarray*}
E(f)=2\sum_{z\in Z}RQ(f_{|lien(x)})ED(f,x).
\end{eqnarray*}
En particulier, si pour tout $z\in Z$ et tout $y\in Y$, $\lambda(lien(z),T_y Y)>1/2$, toute application harmonique é\-qui\-va\-rian\-te $X\to Y$ est constante.
\end{theo}

\preuve
On convient de poser $m(z,z')=0$ si $z$ et $z'$ ne sont pas les extrémités d'une même arête. De même, $m(z,z',z'')=0$ si $z$, $z'$ et $z''$ ne sont pas les sommets d'une même face. On continue de noter $(x',x'')$ (resp. $(x,x',x'')$) un relèvement quelconque à $X$ d'une arête $(z,z')$ (resp. d'une face $(z,z',z'')$). En utilisant la propriété de poids, il vient
\begin{eqnarray*}
E(f)
&=&\frac{1}{2}\sum_{z',z''}m(z',z'')|f(x')-f(x'')|^2 \\
&=&\frac{1}{2}\sum_{z,z',z''}m(z,z',z'')|f(x')-f(x'')|^2 \\
&=&\sum_{z}E(f_{|lien(x)}).
\end{eqnarray*}
D'autre part, pour chaque $x\in X$,
\begin{eqnarray*}
E(f_{|lien(x)})
&=&RQ(f_{|lien(x)})d(f_{|lien(x)},bar(f_{|lien(x)}))^2 \\
&=&RQ(f_{|lien(x)})d(f_{|lien(x)},f(x))^2\\
&=&2RQ(f_{|lien(x)})ED(f,x).
\end{eqnarray*}
Supposons que pour tout $z\in Z$ et tout $y\in Y$, $\lambda(lien(z),T_y Y)\geq\lambda>0$. Alors pour tout $x$, $RQ(f_{|lien(x)})\geq\lambda$, donc
\begin{eqnarray*}
E(f)\geq 2\lambda \sum_{z}ED(f,x)=2\lambda E(f).
\end{eqnarray*}
Si $\lambda>1/2$, cela entraîne que $E(f)=0$, i.e. que $f$ est constante.\qed

\begin{theo}
\label{fixe}
\emph{(H. Izeki et S. Nayatani, \cite{Izeki-Nayatani}).} Soit $Z$ un complexe simplicial fini de groupe fondamental $\Gamma$. Soit $X$ son revêtement universel. Soit $Y$ un espace métrique gé\-o\-dé\-si\-que, $CAT(0)$ et complet. On suppose que pour tout $z\in Z$, 
$$\lambda(lien(z),Y)>\frac{1}{2}.$$
Alors toute action isométrique de $\Gamma$ sur $Y$ possède un point fixe.
\end{theo}

\preuve
Par semi-continuïté du bas du spectre par ultralimite (Lemme \ref{baslim}), l'hypothèse sur les bas du spectre est aussi satisfaite par les cônes asymptotiques $Y'$ de $Y$. La formule de Garland (Théorème \ref{garland}) entraîne que les applications harmoniques é\-qui\-va\-rian\-tes à valeurs dans $Y$ ou l'un de ses cônes asymptotiques sont constantes. On est donc dans le premier cas du Théorème \ref{exharmlim} : il existe un point fixe pour l'action de $\Gamma$ sur $Y$.\qed

\subsection{Propriété (T) de Kazhdan}

\begin{defi}
\label{T}
\emph{(D. Kazhdan, \cite{Kazhdan}).} Un groupe discret $\Gamma$ possède la \emph{propriété (T)} si, pour toute re\-pré\-sen\-ta\-tion unitaire $\rho$ de $\Gamma$ sans vecteurs invariants et tout système générateur fini $S$ de $\Gamma$, il existe une constante $\epsilon(S,\rho)>0$ telle que pour tout vecteur unitaire $\xi\in\hr$,
\begin{eqnarray*}
\max_{s\in S}|\rho(g)\xi -\xi|\geq \epsilon.
\end{eqnarray*}
\end{defi}

La propriété (T) est principalement un phénomène lié à la dimension infinie. On renvoie à \cite{de_la_Harpe-Valette} pour une revue de la propriété (T), et notamment des exemples (non triviaux) suivants.

\begin{exemple}
\label{exT}
Les groupes libres, les groupes moyennables n'ont pas la propriété (T). Les réseaux irréductibles des groupes de Lie semi-simples autres que $PO(n,1)$ et $PU(n,1)$ ont la propriété (T). Les réseaux des groupes $PO(n,1)$ et $PU(n,1)$ ne l'ont pas.
\end{exemple}

\begin{prop}
\label{FH=>T}
Un groupe discret $\Gamma$ a la propriété (T) si et seulement si toute action isométrique affine de $\Gamma$ sur un espace de Hilbert possède un point fixe.
\end{prop}

\preuve
Si $\Gamma$ a un point fixe dans tout espace de Hilbert affine, alors, d'après le lemme \ref{rigcohom}, $H^1 (\Gamma,\rho)=0$, donc l'image du cobord $d:C^{0}(\Gamma,\rho)\to C^{1}(\Gamma,\rho)$ est l'espace des 1-cocycles $Z^1 (\Gamma,\rho)$. Si $\hr$ n'a pas de vecteurs invariants, $d$ est injectif. Soit $S$ un système générateur fini de $\Gamma$. On munit les espaces de cochaînes $C^k (\Gamma,\rho)$ d'une structure hilbertienne en posant 
\begin{eqnarray*}
|\eta|^2 =\sum_{(s_0 ,\ldots,s_k )\in S^{k}}|\eta(s_0 ,\ldots,s_k )|^2 .
\end{eqnarray*}
Alors $d$ est une bijection continue de $C^{0}(\Gamma,\rho)$ sur $Z^{1}(\Gamma,\rho)$. Son inverse est donc continue. On note $N$ la norme de l'inverse. On voit un vecteur unitaire $\xi\in\hr$ comme une 0-cochaîne et $g\mapsto \rho(g)\xi -\xi$ comme le 1-cocycle $d\xi$. Il vient
\begin{eqnarray*}
\max_{s\in S}|\rho(g)\xi -\xi|\geq \frac{1}{\sqrt{|S|}}|d\xi|\geq\frac{1}{N\sqrt{|S|}}|\xi|=\frac{1}{N\sqrt{|S|}}.
\end{eqnarray*}
Par conséquent, $\Gamma$ a la propriété (T).

Pour la réciproque, on renvoie à \cite{de_la_Harpe-Valette}.\qed

\begin{cor}
\label{t}
\emph{(A. Zuk, \cite{Zuk-96}).} Soit $Z$ un complexe simplicial fini de groupe fondamental $\Gamma$. Soit $X$ son revêtement universel. On suppose que pour tout $z\in Z$, $\lambda(lien(z),\R)>1/2$. Alors $\Gamma$ possède la propriété (T).
\end{cor}

\preuve
On applique le Théorème \ref{fixe} aux espaces de Hilbert.\qed

Dans \cite{BS}, W. Ballmann et J. Swiatkowski ont construit des exemples de polyèdres auxquels le corollaire \ref{t} s'applique. On trouvera d'autres exemples dans \cite{Barre} et \cite{DJ}.

\subsection{Inégalité de Garland}

On étend l'identité du théorème \ref{garland}, valable seulement pour les applications harmoniques, en une inégalité valable pour toutes les applications é\-qui\-va\-rian\-tes. Cela permet de prouver directement l'existence d'un point fixe, sans recours aux cônes tangents, et constitue une version quantitative de l'existence d'un point fixe.

On commence par le cas linéaire.

\begin{exemple}
\label{garlandlineaire}
{\em Inégalité de Garland, cas linéaire.} On suppose que, pour tout point $z\in Z$, $\lambda(lien(z),\R)\geq \lambda$. Alors, pour tout application é\-qui\-va\-rian\-te $f:X\to Y$ où $Y$ est un espace de Hilbert, on a l'inégalité
\begin{eqnarray*}
(2\lambda -1)E(f)\leq\lambda \nn\Delta f\nn^{2},
\end{eqnarray*}
où $\nn\Delta f\nn^{2}=\sum_{z}m(z)|\Delta f(x)|^{2}$ est la norme $\ell^{2}$ de $\Delta f$.
\end{exemple}
En effet, on compare deux expressions de l'énergie,
\begin{eqnarray*}
E(f)=\sum_{z}E(f_{|lien(x)})
\end{eqnarray*}
et
\begin{eqnarray*}
E(f)=\sum_{z}ED(f,x)=\sum_{z}\frac{1}{2}\nn f_{|lien(x)}-f(x)\nn^{2}.
\end{eqnarray*}
L'hypothèse sur le bas du spectre des liens donne
\begin{eqnarray*}
E(f)&=&\sum_{z} RQ(f_{|lien(x)})\nn f_{|lien(x)}-bar(f_{|lien(x)})\nn^{2}\\
&\geq& \lambda\sum_{z} \nn f_{|lien(x)}-bar(f_{|lien(x)})\nn^{2}.
\end{eqnarray*}
Pour chaque $x$, la fonction $f_{|lien(x)}-bar(f_{|lien(x)})$ sur $lien(x)$ est orthogonale aux constantes, en particulier à $\Delta f (x)=f(x)-bar(f_{|lien(x)})$, donc, dans $\ell^{2}(lien(x),m)$,
\begin{eqnarray*}
\nn f_{|lien(x)}-f(x)\nn^{2}&=&\nn f_{|lien(x)}-bar(f_{|lien(x)})\nn^{2}+\nn f(x)-bar(f_{|lien(x)})\nn^{2}\\
&=&\nn f_{|lien(x)}-bar(f_{|lien(x)})\nn^{2}+m(x)|\Delta f(x)|^{2}.
\end{eqnarray*}
Il vient
\begin{eqnarray*}
E(f)&\geq& \lambda \sum_{z}\nn f_{|lien(x)}-f(x)\nn^{2}-\lambda \sum_{z}m(x)|\Delta f(x)|^{2}\\
&=&2\lambda E(f)-\lambda \nn\Delta f\nn^{2}.\qed
\end{eqnarray*}

\medskip

Il s'agit de généraliser cette inégalité au cas non linéaire. 

\begin{prop}
\label{ingarland}
Soit $Z$ un complexe simplicial fini de groupe fondamental $\Gamma$. Soit $X$ son revêtement universel. Soit $Y$ un espace métrique géodésique, $CAT(0)$ et complet. Soit $f:X\to Y$ une application é\-qui\-va\-rian\-te. On suppose que, pour tout $z\in Z$ et tout $y\in Y$, $\lambda(lien(z),T_y Y)\geq\lambda>1/2$. Alors
\begin{eqnarray*}
(2\lambda -1)^{2}E(f)\leq 8\lambda^{2}|-\Delta f|^{2}.
\end{eqnarray*}
\end{prop}

\preuve
Reprenons la preuve de la formule de Garland. La première étape est valable pour toute application é\-qui\-va\-rian\-te $f$, 
\begin{eqnarray*}
E(f)=\sum_{z}E(f_{|lien(x)}).
\end{eqnarray*} 
Fixons $x\in X$. Comme la projection $\pi:Y\to T_{f(x)}Y$ diminue les distances et donc les énergies,
\begin{eqnarray*}
E(f_{|lien(x)})\geq E(f_{x}).
\end{eqnarray*}
Par définition du quotient de Rayleigh dans le cône tangent,
\begin{eqnarray*}
E(f_{x})= RQ(f_{x})d(f_x ,bar(f_x ))^2 .
\end{eqnarray*}
Avec l'hypothèse sur le bas du spectre, on a donc
\begin{eqnarray}
\label{minEfx}
E(f)\geq \lambda \sum_{z}d(f_x ,bar(f_x ))^2 .
\end{eqnarray}
Or on se rappelle que
\begin{eqnarray*}
E(f)=\sum_{z}ED(f,x)=\sum_{z}\frac{1}{2}d(f_{|lien(x)},f(x))^{2}.
\end{eqnarray*}
Etant donné $x\in X$, on va donc comparer $d(f_x ,bar(f_x ))^2$ à $2ED(f,x)$. Notons $b_{x}$ (resp. $o_{x}$) l'application constante qui envoie $lien(x)$ sur $bar(f_{x})$ (resp. sur le sommet du cône $T_{f(x)}Y$). Par définition du laplacien, $d(b_{x},o_{x})^{2}=\sum_{x'\in lien(x)}m(x,x')|-\Delta f(x)|^{2}=m(x)|-\Delta f(x)|^{2}$. Comme $\pi$ est radialement isométrique, $d(f_{x},o_{x})^{2}=d(f_{|lien(x)},f(x))^{2}=2ED(f,x)$. Alors
\begin{eqnarray*}
|d(f_x ,bar(f_{x}))^2 -2ED(f,x)|
&=&|d(f_x ,b_{x})^2 -d(f_{x},o_{x})^2 |\\
&=&|d(f_{x},b_{x})+d(f_{x},o_{x})||d(f_{x},b_{x})-d(f_{x},o_{x})|\\
&\leq&2\,d(f_{x},o_{x})d(b_x ,o_{x})\\
&=&2\sqrt{2ED(f,x)}\sqrt{m(x)}|-\Delta f(x)|,
\end{eqnarray*}
par la propriété de minimum du barycentre et l'inégalité triangulaire. Il vient
\begin{eqnarray*}
|\sum_{z}d(f_x ,bar(f_{x}))^2 -2E(f)|
&\leq&\sum_{z} |d(f_x ,bar(f_{x}))^2 -2ED(f,x)|\\
&\leq&\sum_{z} 2\sqrt{2ED(f,x)}\sqrt{m(x)}|-\Delta f(x)|\\
&\leq&2(\sum_{z} 2ED(f,x))^{1/2}(\sum_{z} m(x)|-\Delta f(x)|^{2})^{1/2}\\
&=& 2\sqrt{2E(f)}|-\Delta f|,
\end{eqnarray*}
d'où
\begin{eqnarray*}
\sum_{z}d(f_x ,bar(f_{x}))^2 \geq 2E(f)-2\sqrt{2E(f)}|-\Delta f|.
\end{eqnarray*}
On revient à l'inégalité \ref{minEfx}. Elle donne
\begin{eqnarray*}
E(f)\geq\lambda(2E(f)-2\sqrt{2E(f)}|-\Delta f|),
\end{eqnarray*}
d'où
\begin{eqnarray*}
(2\lambda -1)E(f)-2\lambda\sqrt{2E(f)}|-\Delta f|\leq 0,
\end{eqnarray*}
soit, si $2\lambda -1\geq 0$,
\begin{eqnarray*}
(2\lambda -1)^{2}E(f)\leq 8\lambda^{2}|-\Delta f|^{2}.\qed
\end{eqnarray*}

\begin{theo}
\label{conefixe}
\emph{(H. Izeki et S. Nayatani, \cite{Izeki-Nayatani}).} Soit $Z$ un complexe simplicial fini de groupe fondamental $\Gamma$. Soit $X$ son revêtement universel. Soit $Y$ un espace métrique gé\-o\-dé\-si\-que, $CAT(0)$ et complet. On suppose qu'il existe un $\lambda>\frac{1}{2}$ tel que pour tout $z\in Z$ et tout $y\in Y$, 
$$\lambda(lien(z),T_y Y)\geq\lambda.$$
Alors toute action isométrique de $\Gamma$ sur $Y$ possède un point fixe.
\end{theo}

\preuve
L'inégalité de Garland entraîne qu'on est dans le premier cas du thé\-o\-rè\-me \ref{exharmlim} : il existe un point fixe pour l'action de $\Gamma$ sur $Y$.\qed

\begin{cor}
\label{lissefixe}
\emph{(Compare M.T. Wang, \cite{Wang1}).} Soit $Z$ un complexe simplicial fini de groupe fondamental $\Gamma$. Soit $X$ son revêtement universel. Soit $Y$ une variété riemannienne complète, simplement connexe, à courbure sectionnelle négative ou nulle. On suppose que pour tout $z\in Z$, $\lambda(lien(z),\R)>1/2$. Alors toute action isométrique de $\Gamma$ sur $Y$ possède un point fixe.
\end{cor}

\section{Calculs de bas de spectre}
\label{basspectre}

Dans cette section, on donne des exemples de bas de spectre pour des graphes munis du poids naturel : toutes les arêtes ont un poids égal à 1, et les sommets un poids égal à leur valence. On commence par rappeler des résultats classiques sur les bas de spectres scalaires, puis on montre comment les cas des variétés lisses puis des arbres se ramènent au cas scalaire. Enfin, on indique deux tentatives infructueuses pour minorer des bas de spectres à valeurs dans des immeubles.

\subsection{Bas du spectre scalaire des immeubles}

\begin{exemple}
\label{bascycle}
Si $C_k$ est un cycle de longueur $k$, 
\begin{eqnarray*}
\lambda(C_k ,\R)=\frac{1}{2}|1-e^{2i\pi/k}|^2 .
\end{eqnarray*}
\end{exemple}
En effet, on utilise la transformation de Fourier discrète pour diagonaliser la forme quadratique énergie. Les valeurs propres sont les $\frac{1}{2}|1-\zeta|^2$ où $\zeta$ parcourt les racines $k$-èmes de l'unité, et $\zeta=1$ correspond aux fonctions constantes.

On ne trouve $\lambda(C_k ,\R)>1/2$ que si $k\leq 5$. Un complexe simplicial dont les liens sont des cycles de longueur $\leq 5$ est une surface à courbure positive, dont le groupe fondamental est fini. La propriété de point fixe pour tout espace $CAT(0)$ est immédiate, le théorème \ref{fixe} n'apporte pas grand chose. On trouve que $\lambda(C_6 ,\R)=1/2$, c'est le cas limite pour l'application du théorème \ref{fixe}. On essaiera d'en tirer parti au paragraphe \ref{sspec}.

\begin{defi}
\label{defplan}
On appelle \emph{triangle généralisé} un graphe ayant la propriété suivante~: tout cycle est de longueur au moins 6, et deux arêtes quelconques sont contenues dans un cycle de longueur 6. 
\end{defi}

\begin{exemple}
\label{explan}
Soit $\mathbf{F}$ un corps fini à $q$ éléments. Considérons le graphe biparti dont l'ensemble des sommets est la réunion de l'ensemble des points et de l'ensemble des droites du plan projectif $\mathbf{F}P^2$. On met une arête entre un point $p$ et une droite $d$ si $p\in d$. Il s'agit d'un triangle généralisé de valence $q+1$.
\end{exemple}
En effet, un cycle correspond à un polygone, qui a au moins 3 côtés, donc le cycle est de longueur au moins 6. Deux paires $(p\in d)$ et $(p'\in d')$ sont toujours contenues dans le triangle dont le troisième sommet est $p''=d\cap d'$.

\begin{prop}
\label{lambdatri}
\emph{(W. Feit, G. Higman \cite{Feit-Higman}).} Si $C$ est un triangle généralisé de valence $q$, alors
\begin{eqnarray*}
\lambda(C,\R)=1-\frac{\sqrt{q-2}}{q-1}.
\end{eqnarray*}
\end{prop}

Les {\em immeubles de type} $\tilde{A}_{2}$ sont les complexes simpliciaux dont les liens sont des triangles généralisés. Le théorème \ref{fixe} et la propriété \ref{FH=>T} entraînent que les groupes discrets cocompacts d'automorphismes des immeubles de type $\tilde{A}_{2}$ ont la propriété de point fixe sur les espaces de Hilbert, i.e. la propriété (T), \cite{Pansu1}, \cite{Pansu2}, \cite{Zuk-96}. C'était déjà connu par d'autres méthodes pour beaucoup d'entre eux, mais pas tous, voir \cite{CMSZ}.

\subsection{Bas du spectre scalaire d'un graphe générique}

\emph{``La plupart des graphes ont un bas du spectre proche de 1, lorsque le nombre de sommets est grand par rapport à la valence''}. Comment préciser cet énoncé ? 

Tout graphe à $n$ sommets de valence $2d$ peut-être obtenu de la façon suivante~: étant donné un ensemble $S$ à $n$ éléments, on se donne $d$ permutations $\sigma_1 ,\ldots,\sigma_d$ de $S$, on construit un graphe $\gamma(\sigma_1 ,\ldots,\sigma_d )$ de valence $2d$ dont l'ensemble des sommets est $S$ en reliant chaque sommet $s\in S$ aux sommets $\sigma_1 (s),\ldots,\sigma_d (s)$. Noter qu'on accepte les arêtes multiples (cas où $\sigma_i (s)=\sigma_j(s)$ ou $\sigma_{j}^{-1}(s)$) et les boucles reliant un sommet à lui-même (cas où $\sigma_i (s)=s$).

\def\sn{\mathcal{S}_{n}}
On utilise la mesure de probabilité unifome sur les $d$-uplets de permutations (le produit de $d$ copies du groupe symétrique $\sn$), et la mesure image par l'application $\gamma$ sur les graphes. On obtient ainsi une loi de probabilité sur l'ensemble $L(n,d)$ des graphes à $n$ sommets de valence $2d$. On choisit le poids $m$ uniforme qui vaut 1 sur chaque arête et donc $2d$ en chaque sommet. 

\begin{theo}
\label{friedman}
{\em (A. Broder, E. Shamir, \cite{Broder-Shamir}).} Lorsque $n$ tend vers l'infini, la proportion de graphes dans $L(n,d)$ dont le bas du spectre est proche de 1 tend vers 1. Plus précisément {\em (J. Friedman, \cite{Friedman})}, il existe une constante $c$ indépendante de $k$ telle que la proportion de graphes $C$ dans $L(n,d)$ dont le bas du spectre satisfait
\begin{eqnarray*}
\lambda(C,\R)\geq 1-(\frac{\sqrt{2d-1}}{d}+\frac{\log(2d)}{2d}+\frac{c}{d})
\end{eqnarray*}
tend vers 1 quand $n$ tend vers l'infini.
\end{theo}

\preuve
On donne seulement les grandes lignes de la preuve du résultat de Broder et Shamir.

1. Concentration. La fonction $\lambda:(\sn)^d \to\R$ qui à un $d$-uplet de permutations $\sigma_1 ,\ldots,\sigma_d $, associe le bas du spectre du graphe $\gamma(\sigma_1 ,\ldots,\sigma_d )$, est lipschitzienne au sens suivant. Si deux $d$-uplets $\mathbf{\sigma}$ et $\mathbf{\sigma'}$ ne diffèrent que par une des permutations, alors
\begin{eqnarray*}
|\lambda(\mathbf{\sigma})-\lambda(\mathbf{\sigma'})|\leq \frac{4}{d}.
\end{eqnarray*}
En effet, soit $g:S\to\R$ une fonction telle que $bar(g)=0$. Le dénominateur $\displaystyle d(g,bar(g))^2 =2d\sum_s |g(s)-bar(g)|^2$ du quotient de Rayleigh ne dépend pas des arêtes. Quant au numérateur, il vaut
\begin{eqnarray*}
E(g)&=&\frac{1}{2}\sum_{s\in S}\sum_{i=1}^d (|g(s)-g(\sigma_i (s))|^2 +|g(s)-g(\sigma_{i}^{-1} (s))|^2 )\\
&=&\sum_{s\in S}\sum_{i=1}^d |g(s)-g(\sigma_i (s))|^2 .
\end{eqnarray*}
Si deux $d$-uplets $\mathbf{\sigma}$ et $\mathbf{\sigma'}$ ne diffèrent que par la $j$-ème permutation, alors
\begin{eqnarray*}
|E_\mathbf{\sigma}(g)-E_\mathbf{\sigma'}(g)|
&\leq &\frac{1}{2}\sum_{s\in S}(|g(s)-g(\sigma_j (s))|^2 -|g(s)-g({\sigma'}_{j}^{-1} (s))|^2 )\\
&\leq &\sum_{s\in S}2(|g(s)|^2 +|g(\sigma_j (s))|^2 +|g(s)|^2 +|g(\sigma'_j (s))|^2 )\\
&\leq &8\sum_{s\in S}|g(s)|^2 \\
&=&\frac{4}{d}d(g,bar(g))^2 .
\end{eqnarray*}
En prenant la borne inférieure sur toutes les fonctions $g$ de moyenne nulle, on obtient l'inégalité annoncée, laquelle énonce que la fonction $\lambda$ est 1-lipschitzienne pour la distance sur $(\sn)^d$ produit (au sens $\ell^{1}$) de la distance sur $\sn$ qui met deux permutations quelconques à distance $4/d$ l'une de l'autre.

D'un résultat de concentration (voir \cite{Ledoux}, Corollary 1.17), il résulte que, hors d'un sous-ensemble de $(\sn)^d$ de mesure exponentiellement petite, la fonction $\lambda$ est très proche de son espérance. Précisément, pour tout $r>0$,
\def\E{\mathbf{E}}
\begin{eqnarray*}
P(\lambda<\E(\lambda)-r)<e^{-r^2 /2D^{2}},
\end{eqnarray*}
où $D^{2}$ est la somme des carrés des diamètres de chaque copie de $(\sn)^d$, soit $d\times (4/d)^{2}=16/d$. En prenant $r$ de l'ordre de $\frac{1}{\sqrt{d}}$, on voit qu'il suffit de minorer l'espérance de la variable $\lambda$.

2. Lien entre bas du spectre et probabilité de retour. On considère la marche aléatoire simple sur le graphe $\gamma(\sigma)$. La matrice $P$ des probabilités de transition est $P(s,s')=\frac{1}{k}$ s'il existe $i$ tel que $s'=\sigma_{i}(s)$ ou $s'=\sigma_{i}^{-1}(s)$, $P(s,s')=0$ sinon. La probabilité qu'une marche aléatoire issue de $s$ revienne en $s$ au bout de $2k$ pas est le coefficient de matrice $P^{2k} (s,s)$. La probabilité qu'une marche aléatoire partie d'un point quelconque y revienne au bout de $2k$ pas est égale à la trace $tr(P^{2k})$. Notons $1=\mu_1 \geq \cdots\geq\mu_n$ les valeurs propres de $P$. Alors le bas du spectre $\lambda=1-\mu_2$, d'où
\begin{eqnarray*}
tr(P^{2k})=1+\mu_2^{2k} +\cdots+\mu_n^{2k} \geq 1+\mu_2^{2k} =1+(1-\lambda)^{2k}.
\end{eqnarray*}
On voit $\lambda$ comme une variable aléatoire sur $(\sn)^d$. Il vient
\begin{eqnarray*}
1-\E(\lambda)=\E(1-\lambda)\leq \E((1-\lambda )^{2k})^{1/2k}\leq (tr(P^{2k})-1)^{1/2k}.
\end{eqnarray*}

\def\M{\mathcal{M}}
3. Estimation de l'espérance de la probabilité de retour. Faire un pas au hasard à partir de $s$ consiste à tirer au hasard l'un des points $\sigma_{i}^{\pm}(s)$. Marcher au hasard pendant le temps $2k$ à partir de $s$ consiste à choisir indépendamment $2k$ indices dans $\{1,\ldots,d\}\times\pm 1$. Un tel choix peut-être vu comme une marche dans le monoïde libre $\M$ engendré par $2d$ symboles $\{\pi_1 ,\pi_1^{-1},\ldots,\pi_d ,\pi_d^{-1}\}$. Marchons maintenant au hasard dans un graphe tiré au hasard. On intervertit les tirages aléatoires. On tire d'abord le mot $\mathbf{\pi}$ de longueur $2k$ dans $\M$, puis le $d$-uplet $\mathbf{\sigma}$ dans $(\sn)^d$. Fixons un sommet $s_0$. Alors $\displaystyle \frac{1}{n}(n!)^d (2d)^{2k}\E(tr(P^{2k}))$ est égal au nombre de tirages $\mathbf{\pi}$ et $\mathbf{\sigma}$ donnant une marche qui part et se termine en $s_0$.

Si le mot $\mathbf{\pi}$ n'est pas réduit dans le groupe libre de générateurs $\{\pi_1 ,\ldots,\pi_d \}$, la marche fait des aller-retour sans intérêt. On majore aisément le nombre de mots non réduits. La contribution principale vient des mots triviaux dans le groupe libre. Leur probabilité est contrôlée par le bas du spectre d'un arbre régulier de valence $2d$, elle est inférieure à $(2/d)^k$.

Désormais, on suppose le mot $\mathbf{\pi}$ réduit. La contribution principale va venir des marches qui ne passent jamais deux fois par le même sommet jusqu'au temps $2k-1$. En effet, lorsque, à $\pi$ fixé, on compte le nombre de $d$-uplets $\mathbf{\sigma}$ qui réalisent une trajectoire donnée, chaque sommet rencontré pour la première fois donne une seule contrainte sur une permutation $\sigma_i$, alors qu'un sommet déjà rencontré (auto-intersection) en impose deux. La probabilité qu'une marche se recoupe avant de revenir en $s_0$ est donc au moins $n$ fois plus petite que la probabilité qu'une marche revienne en $s_0$ sans se recouper, laquelle est inférieure à $\displaystyle \frac{1}{n-k}$.

On obtient donc une majoration de la forme $\displaystyle\E(tr(P^{2k}))\leq n((\frac{2}{d})^k +\frac{1}{n-k}+reste)$. On choisit $k$ de l'ordre de $\log_d (n)$ et on conclut que $\E(\lambda)$ se comporte asymptotiquement en $d$ comme pour l'arbre régulier de valence $2d$, pour lequel $\displaystyle \lambda=1-\frac{2\sqrt{2d-1}}{2d}$.\qed

\begin{rem}
Le modèle avec les permutations peut paraître artificiel. Il est simplement commode. Une modification légère permet d'étendre la conclusion du thé\-o\-rè\-me \ref{friedman} à la mesure uniforme sur l'ensemble $L(n,d)$. On peut même éliminer les graphes possédant des arêtes multiples ou des boucles, voir \cite{Broder-Shamir}.
\end{rem}

\subsection{Cônes tangents}

Dans l'espace euclidien, si $s(t)=y+tv$ et $s'(t)=y+tv'$ sont deux demi-droites issues de $y$, on a pour tout $t$ l'identité
\begin{eqnarray*}
|v-v'|=\frac{d(s(t),s'(t))}{t}.
\end{eqnarray*}

Dans un espace $CAT(0)$, étant données deux géodésiques $s$ et $s'$ d'origine $y$, parcourues à vitesse constante, la fonction $t\mapsto \frac{d(s(t),s'(t))}{t}$ est croissante. On peut donc poser
\begin{eqnarray*}
d(s,s')=\lim_{t\to 0} \frac{d(s(t),s'(t))}{t},
\end{eqnarray*}
et identifier $s$ et $s'$ si $d(s,s')=0$. 

\begin{defi}
\label{defcone}
\emph{(\cite{Bridson-Haefliger})}.
On obtient un espace métrique appelé {\em cône tangent} de $Y$ en $y$ et noté $T_{y}Y$. 
\end{defi}

C'est à nouveau un espace $CAT(0)$, car c'est une limite d'espaces $CAT(0)$. Il vient avec une application $\pi:Y\to T_{y}Y$ qui est isométrique le long de chaque géodésique issue de $y$, et diminue les distances en général. $\pi$ envoie un point $y'$ sur le segment géodésique de $y$ à $y'$ paramétré à vitesse constante sur $[0,1]$.

\begin{exemple}
\label{excone}
Si $Y$ est une variété riemannienne simplement connexe à courbure négative ou nulle (par exemple, un espace symétrique de type non compact), ses cônes tangents sont des espaces euclidiens. 

Si $Y$ est un arbre, ses cônes tangents sont des réunions de demi-droites attachées à leur extrémité.
\end{exemple}

\begin{prop}
\label{bastan}
{\em (Wang,\cite{Wang2}).} Pour tout complexe simplicial pondéré fini $(C,m)$ et tout espace $CAT(0)$ $Y$,
\begin{eqnarray*}
\lambda(C,Y) = \inf_{y\in Y} \lambda(C,T_{y}Y).
\end{eqnarray*}
\end{prop}

\preuve
Etant donnée $g:C\to Y$, on note $g'=\pi\circ g$ où $\pi$ est la projection $Y\to T_{bar(g)}Y$. Alors $bar(g')$ est le sommet du cône. En effet, si $t\mapsto s(t)$ est une géodésique dans $Y$, issue de $bar(g)$ et parcourue à vitesse constante, et si $y\in Y$, alors la dérivée en $t=0$ de la fonction $t\mapsto d(y,s(t))$ s'exprime en fonction des points $\pi(s(1))$ et $\pi(y)$ du cône tangent  (\cite{Bridson-Haefliger}, corollaire II.3.6). Comme le cône tangent au cône en son sommet est le cône lui-même, $\pi\circ\pi=\pi$, et les fonctions $t\mapsto d(y,s(t))$ et $t\mapsto d(\pi(y),\pi(s(t)))$ ont même dérivée en $t=0$. En faisant la moyenne sur $C$, on trouve que les fonctions $t\mapsto d(g,s(t))^2$ et $t\mapsto d(\pi\circ g,\pi(s(t)))^2$ ont même dérivée en $t=0$. Comme la première atteint un minimum en 0, il vient
\begin{eqnarray*}
\frac{d}{dt}{d(\pi\circ g,\pi(s(t)))^2 }_{|t=0}\geq 0.
\end{eqnarray*}
Comme cette fonction est convexe, elle atteint aussi un minimum en 0. Cela prouve que $y\mapsto d(g',\pi(y))^2$ atteint son minimum en $bar(g)$, donc $bar(g')$ est le sommet du cône.
 
Comme $\pi$ diminue les distances, $E(g')\leq E(g)$. Comme $\pi$ est isométrique le long des géodésiques issues de $bar(g)$, $d(g',\pi(bar(g)))^2=d(g,bar(g))^2$. Par conséquent, $RQ(g')\leq RQ(g)$. Cela prouve que $\lambda(C,Y) \geq \inf_{y\in T}\lambda(C,T_{y}Y)$.

L'inégalité inverse résulte de fait que les cônes tangents sont des ultralimites particulières, et du Lemme \ref{baslim}.\qed

\begin{exemple}
\label{exbastan}
Si $Y$ est une variété riemannienne simplement connexe à courbure négative ou nulle (e.g. un espace symétrique sans facteurs compacts), alors 
\begin{eqnarray*}
\lambda(C,m,Y)=\lambda(C,m,\R).
\end{eqnarray*}
\end{exemple}

C'est en utilisant cette propriété du bas du spectre que M.T. Wang prouvait le Théorème \ref{conefixe} et le Corollaire \ref{lissefixe} dans le cas des actions réductives.

\subsection{Ultralimites}

La semi-continuïté du bas du spectre par ultralimite (Lemme \ref{baslim}) permet d'obtenir quelques valeurs supplémentaires.

\begin{cor}
\label{bascone}
\emph{(Gromov, \cite{Gromov-RWRG}).}
Soit $M$ une variété riemannienne simplement con\-nexe, à courbure sectionnelle négative ou nulle. Soit $cone_{\omega}(M)$ un cône asymptotique de $M$. Alors pour tout graphe pondéré fini $(C,m)$,
\begin{eqnarray*}
\lambda(C,cone_{\omega}(M))=\lambda(C,\mathbf{R}).
\end{eqnarray*}
\end{cor}

\begin{exemple}
\label{basconesym}
\emph{(B. Kleiner, B. Leeb, \cite{Kleiner-Leeb})}.
Les cônes asymptotiques des espaces sy\-mé\-tri\-ques sont des immeubles euclidiens. Pour de tels immeubles $I$, on a donc, pour tout graphe pondéré fini $(C,m)$,
\begin{eqnarray*}
\lambda(C,I)= \lambda(C,\mathbf{R}).
\end{eqnarray*}
\end{exemple}

\begin{rem}
\label{remimm}
Cette égalité ne peut être vérifiée par tous les immeubles euclidiens.
\end{rem}
En effet, soit $I$ un immeuble euclidien localement fini de dimension $d\geq 2$. Soit $C=lien(x)$ le lien d'un sommet, muni de la pondération naturelle, qui donne le poids 1 aux simplexes de dimension $d$. Il est fréquent que $\lambda(C,\mathbf{R})>1/2$ (c'est le cas pour les immeubles de type $\tilde{A}_2$, d'après la proposition \ref{lambdatri} ; pour les autres immeubles, voir par exemple \cite{Garland}). Il est aussi fréquent que $I$ possède un groupe discret cocompact d'isométries. Si on avait $\lambda(C,I)\geq \lambda(C,\mathbf{R})$, on prouverait à l'aide du théorème \ref{garland} que toute application harmonique é\-qui\-va\-rian\-te $I\to I$ est constante. Or l'identité $id:I\to I$ est harmonique non constante, contradiction.

\subsection{Bas du spectre à valeurs dans un arbre}

On rappelle que, d'après la proposition \ref{bastan}, pour minorer $\lambda(C,m,Y)$ pour un espace $CAT(0)$ $Y$, il suffit de traiter ses cônes tangents. Plus précisément, il suffit de minorer $RQ(g)$ lorsque $g:C\to Y$ a son barycentre au sommet du cône.

\begin{prop}
\label{lamarbre}
\emph{(H. Izeki et S. Nayatani, \cite{Izeki-Nayatani}).} Si $Y$ est un arbre, alors pour tout graphe pondéré $(C,m)$,
\begin{eqnarray*}
\lambda(C,Y)=\lambda(C,\R).
\end{eqnarray*}
\end{prop}

\preuve
Comme on n'a affaire qu'à des ensembles finis de points de $Y$, on peut supposer que $Y$ est un arbre fini.

Un cône tangent à un arbre fini est ou bien une droite (dans ce cas, il n'y a rien à prouver), ou bien une réunion de demi-droites $D_1 ,\ldots,D_n$ de même origine. Soit $g:C^0 \to\bigcup_{i}D_{i}$ une application dont le barycentre est au sommet. Notons $a_i$ la somme pondérée des distances au sommet des points de $g(C^0 )$ se trouvant dans la demi-droite $D_i$. Comme $bar(g)$ est au sommet, pour tout $i=1,\ldots,n$,
\begin{eqnarray*}
a_i \leq \sum_{j\not=i} a_j .
\end{eqnarray*}
En effet, lorsqu'on avance dans la branche $D_i$, la dérivée de la fonction $d(g,\cdot)^2$ à l'origine est $-a_i +\sum_{j\not=i} a_j \geq 0$. Par conséquent, il existe un polygone dans le plan euclidien dont les longueurs des côtés sont les $a_i$. Autrement dit, il existe des vecteurs unitaires $e_1 ,\ldots,e_n$ dans $\R^2$ tels que $\sum_{i=1}^{n} a_i e_i =0$. On définit une application $g':C\to\R^2$ en posant $g'(c)=d(g(c),bar(g))e_i$ lorsque $g(c)\in D_i$. Alors $bar(g')=0$, $E(g')\leq E(g)$ et $d(g',bar(g'))^2 =d(g,bar(g))^2$, donc $RQ(g')\leq RQ(g)$. On conclut que $\lambda(C,m,Y)\geq\lambda(C,m,\R^2)=\lambda(C,m,\R)$, puis, comme $Y$ contient des copies isométriques de $\R$, que $\lambda(C,m,Y)=\lambda(C,m,\R)$.\qed

\subsection{L'invariant d'Izeki-Nayatani}

Ce paragraphe développe l'idée illustrée dans le cas des arbres de remplacer une application à valeurs dans $Y$ par une application à valeurs dans un espace euclidien.

\begin{defi}
{\em (H. Izeki and S. Nayatani, \cite{Izeki-Nayatani})}.
Soit $Y$ un espace géodésique $CAT(0)$ complet. Etant donné un sous-ensemble pondéré fini $Z \in Y$ (on suppose que la somme des poids vaut 1), soit $\phi:Z\to\mathcal{H}$ une application 1-Lipschitzienne à valeurs dans un espace de Hilbert telle que pour tout $z\in Z$, $|\phi(z)|=d(z,bar(Z))$. On définit
\begin{eqnarray*}
\delta(Z)=\inf_{\phi}\frac{|bar(\phi)|^2}{\nn\phi\nn^2}.
\end{eqnarray*} 
L'\emph{invariant d'Izeki-Nayatani} de $Y$ is $IN(Y)=\sup_{Z\subset Y}\delta(Z)\in[0,1]$.
\end{defi}

\begin{lemme}
\label{in}
Soit $Y$ un espace géodésique $CAT(0)$ complet, soit $C$ un graphe fini pondéré. Alors
\begin{eqnarray*}
\lambda(C,Y)\geq (1-IN(Y))\lambda(C,\R).
\end{eqnarray*}
\end{lemme}

\preuve
Etant donnée une application $g:C^{0}\to Y$, on pose $Z=g(C)$. On choisit une application optimale $\phi$ for $Z$. D'après Pythagore, 
$$d(\phi,bar(\phi))^2 =\nn\phi\nn^2 -|bar(\phi)|^2=(1-\delta(Z))\nn\phi\nn^2=(1-\delta(Z))d(g,bar(g))^2.$$
Alors
\begin{eqnarray*}
\lambda(C,\R)&\leq& RQ(\phi\circ g)=\frac{E(\phi\circ g)}{d(\phi\circ g,bar(\phi\circ g))^2}\\
&\leq& \frac{E(g)}{d(\phi,bar(\phi))^2}=\frac{1}{1-\delta(Z)}RQ(g).
\end{eqnarray*}

\begin{exemple}
\label{exin}
\begin{enumerate}
  \item Par définition, pour un espace de Hilbert, $IN=0$.
  \item Pour tout $Y$, $IN(Y)=\inf_{y\in Y}IN(T_y Y)$. Par conséquent, les variétés riemanniennes complètes simplement connexes à courbure sectionnelle négative ou nulle satisfont $IN=0$.
\item La preuve de la Proposition \ref{lamarbre} donne $IN=0$ pour les arbres.
\item $IN$ est continu pour les ultralimites. Par conséquent, les immeubles euclidiens qui apparaissent comme cônes asymptotiques d'espaces sy\-mé\-tri\-ques satisfont $IN=0$.
\item Pour tout $Y$ et tout espace probabilisé $\Omega$, $IN(L^2 (\Omega,Y))\leq IN(Y)$.
\item $IN(Y_1 \times Y_2 )\leq \max\{IN(Y_1 ),IN(Y_2 )\}$. Par conséquent, $IN=0$ pour tous les produits dont les facteurs appartiennent à l'une des familles 1 à 5.
\item L'immeuble de Bruhat-Tits $I_p$ associé au groupe $Sl(3,\Q_p )$ satisfait 
$$IN(I_p )\geq\frac{(\sqrt{p}-1)^2}{2(p-\sqrt{p}+1)}.$$
On conjecture que l'égalité a lieu.
\end{enumerate}
\end{exemple}
La minoration de $IN$ pour les immeubles $I_p$ résulte du Lemme \ref{in}, du Lemme \ref{lambdatri}, du théorème \ref{fixe}, et du fait qu'il existe une application harmonique non constante de $I_p$ dans $I_p$, à savoir l'identité.

\begin{theo}
\label{sl3}
\emph{(H. Izeki et S. Nayatani, \cite{Izeki-Nayatani}).} 
Soit $I_2$ l'immeuble de Bruhat-Tits associé au groupe $Sl(3,\Q_2 )$. Alors
\begin{eqnarray*}
IN(I_2 )\leq 0.4122.
\end{eqnarray*}
Par conséquent, pour tout graphe pondéré fini $C$,
$$\lambda(C,I_{2})\geq 0.5878\,\lambda(C,\R).$$
\end{theo}

\preuve
Elle n'est pas simple (on renvoie à \cite{Izeki-Nayatani}) et donne une majoration de $IN(I_p )$ qui est malheureusement $>1/2$ pour $p\geq 3$.\qed

\begin{defi}
\label{deffy}
\emph{(H. Izeki, T. Kondo et S. Nayatani, \cite{Izeki-Kondo-Nayatani}).} 
Soit $\delta \in [0,1]$. On dit qu'un groupe $\Gamma$ a la \emph{propriété} $F\mathcal{Y}_{\leq\delta}$ si toute action isométrique de $\Gamma$ sur un espace géodésique $CAT(0)$ complet $Y$ tel que $IN(Y)\leq \delta$ a un point fixe.
\end{defi}

\begin{prop}
\label{propfy}
Soit $\delta <\frac{1}{2}$. Soit $Z$ un complexe simplicial fini de groupe fondamental $\Gamma$. On suppose que pour tout $z\in Z$, $\lambda(lien(z),\R)>\frac{1}{2(1-\delta)}$. Alors $\Gamma$ a la propriété $F\mathcal{Y}_{\leq\delta}$.
\end{prop}

Les classes $F\mathcal{Y}_{\leq\delta}$ sont intéressantes. En effet, pour la topologie de l'espace des groupes marqués (voir \cite{Champetier}), ces classes sont ouvertes, \cite{Kondo}.

Le minorant de $IN(I_p )$ donné en \ref{exin} tend vers $1/2$ quand $p$ tend vers l'infini. Par conséquent, chaque classe $\mathcal{Y}_{\leq\delta}$, $\delta <\frac{1}{2}$, ne contient qu'un nombre fini de ces immeubles. La Proposition \ref{propfy} est donc loin de conduire à la propriété FSI. Néanmoins, pour tout $\delta <\frac{1}{2}$, la classe $F\mathcal{Y}_{\leq\delta}\cap FSI$ est dense dans l'espace des groupes marqués, voir \cite{Kondo}.

\subsection{Inégalités de Wirtinger}
\label{sspec}

Ce paragraphe et les deux suivants proposent une autre approche, destinée donner un minorant uniforme du bas du spectre à valeurs dans un immeuble. Elle s'appuie sur un résultat de M. Gromov sur des inégalités satisfaites par tous les espaces géodésiques $CAT(0)$, qu'il a baptisées inégalités de Wirtinger.

\begin{defi}
Soit $C$ un graphe fini non pondéré, soit $Y$ un espace mé\-tri\-que, soit $g:C\to Y$. Pour $j=1,\,2,\ldots$, on note
\begin{eqnarray*}
E_{j}(g)=\frac{1}{2}\sum_{c\in C}\sum_{d(c,c')=j} d(g(c),g(c'))^{2}.
\end{eqnarray*}
\end{defi}
Autrement dit, pour $j=1$, $E_1$ est l'énergie relative à la pondération uniforme (chaque arête a un poids égal à 1). 

\begin{exemple}
Cas du $k$-cycle $C_{k}$.
\end{exemple}
Soit $C_{k}$ le graphe formé d'un unique cycle de $k$ arêtes, avec la pondération uniforme. Soit $g:C_{k}^{0}\to\R^{2}$ un plongement de $C_{k}^{0}$ comme un polygone régulier à $k$ côtés dans le plan euclidien, inscrit dans un cercle de rayon 1. On note $W(k,1)=E(g)=4k\sin^{2}(\pi/k)$ son énergie. Alors $E_{j}(g)=W(k,j):=4k\sin^{2}(\pi j/k)$.

\begin{center}
\includegraphics[width=2in]{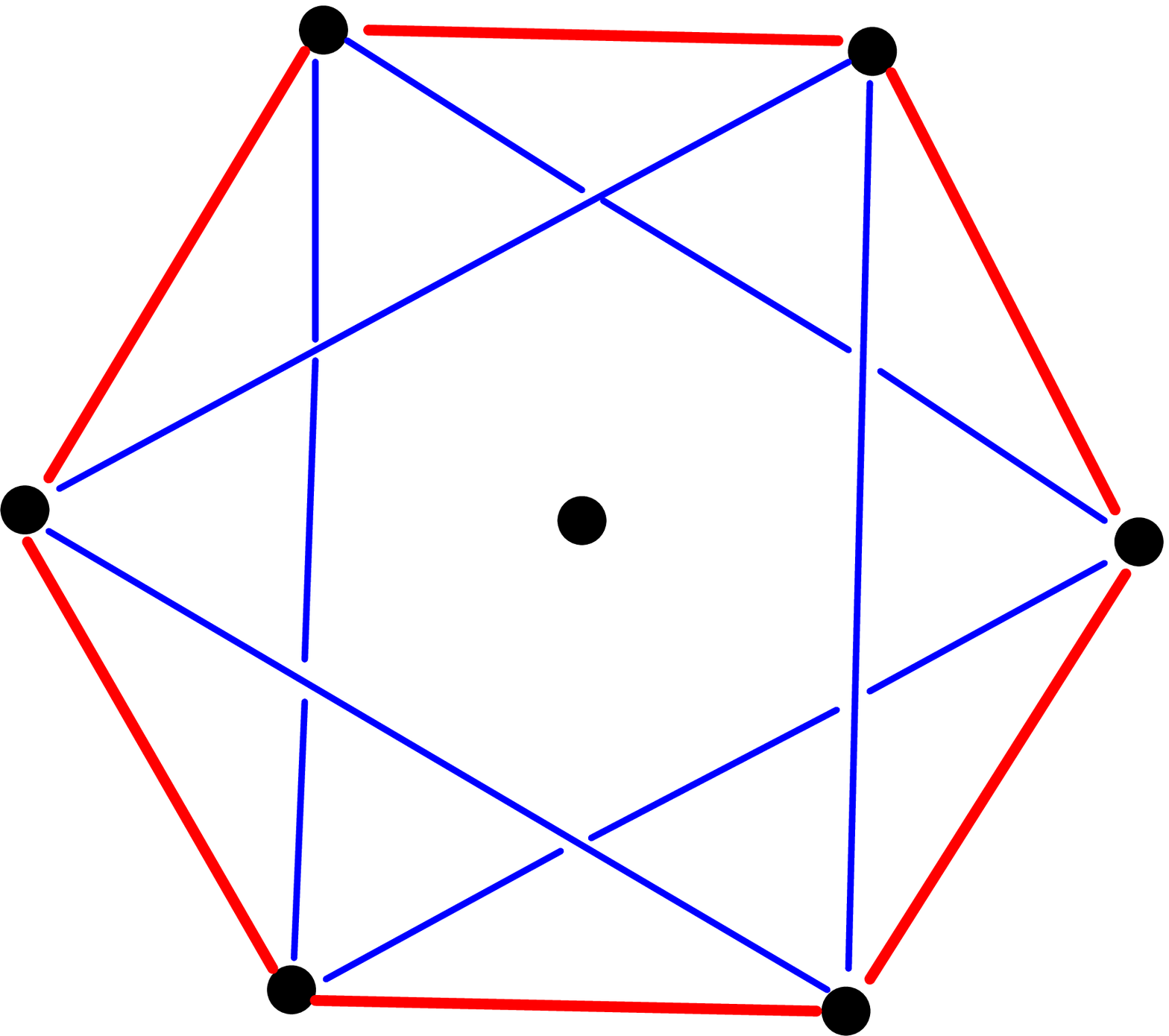}

En rouge, le graphe $C_6$. En bleu, les couples de points entrant dans la définition de $E_{2}$.
\end{center}  

\begin{defi}
Suivant Gromov \cite{Gromov-CAT}, on dit qu'un espace métrique $Y$ satisfait à {\em l'inégalité de Wirtinger} $Wir_k$ si, pour toute application $g:C_{k}\to Y$, et pour tout $j=1,\ldots,k$,
\begin{eqnarray*}
\frac{E(g)}{E_{j}(g)}\geq \frac{W(k,1)}{W(k,j)}.
\end{eqnarray*}
\end{defi}

\begin{prop}
\label{rnwir}
Les espaces euclidiens satisfont à toutes les inégalités de Wirtinger $Wir_{k}$, $k\geq 4$.
\end{prop}

\preuve
Soit $g:C_{k}\to\R^{n}$. On identifie $C_{k}$ à l'ensemble des racines $k$-èmes de l'unité dans $\C$, et on note $\zeta=e^{2i\pi/k}$. On décompose $g$ en série de Fourier discrète $g=\sum_{\ell=1}^{k}g_{\ell}$, où $g_{\ell}(\zeta c)=\zeta^{\ell}g_{\ell}(c)$. D'après l'identité de Bessel-Parseval, cette décomposition est orthogonale pour la norme $\ell^{2}$,
\begin{eqnarray*}
\sum_{c}|g(c)|^{2}=\sum_{\ell}\sum_{c}|g_{\ell}(c)|^{2}.
\end{eqnarray*}
Comme
\begin{eqnarray*}
E_{j}(g)&=&\sum_{c}|g(\zeta^{j}c)-g(c)|^{2}\\
&=&\sum_{c}|\sum_{\ell}(g_{\ell}(\zeta^{j}c)-g_{\ell}(c))|^{2}\\
&=&\sum_{c}|\sum_{\ell}(\zeta^{j\ell}-1)g_{\ell}(c)|^{2}\\
&=&\sum_{\ell}|\zeta^{j\ell}-1|^{2}\sum_{c}|g_{\ell}(c)|^{2},
\end{eqnarray*}
la forme hermitienne $E_{j}$ est diagonale. Par conséquent, la borne inférieure des nombres $E_{1}(g)/E_{j}(g)$ est égale au plus petit des quotients $\displaystyle \frac{|\zeta^{\ell}-1|^{2}}{|\zeta^{j\ell}-1|^{2}}$, $\ell=0,\ldots,k-1$. On vérifie aisément que le quotient est minimum exactement lorsque $\ell=1$ ou $\ell=k-1$. La fonction $g:C_{k}\to\R^{2}$ définie par $g(c)=(\Re e(c),\Im m(c))$ réalise ce minimum, donc celui-ci est égal à $\frac{W(k,1)}{W(k,j)}$.

Pour un usage ultérieur, remarquons que, pour une application $g:C_{k}\to \R^{n}$, $E_{1}(g)/E_{j}(G)=\frac{W(k,1)}{W(k,j)}$ si et seulement si il existe des vecteurs $u\in\R^{n}$ et $v\in\C^{n}$ tels que $g(c)=u+cv+\overline{cv}$, i.e. si $g$ est la restriction aux racines de l'unité d'une application $\R$-affine de $\C$ dans $\R^{n}$.\qed

\begin{theo}
\label{cat=>wir}
{\em (M. Gromov, \cite{Gromov-CAT} section 25)}. Tout espace $CAT(0)$ satisfait aux inégalités de Wirtinger $Wir_{k}$ pour tout $k\geq 4$.
\end{theo}

\preuve
Etant donnée $g:C_{k}\to Y$, on équipe chaque arête $cc'$ de $C_{k}$ de la métrique qui lui donne la longueur $d(g(c),g(c'))$, de sorte que $C_{k}$ devient un cercle $\partial D$. On prolonge $g$ en une application isométrique $f:\partial D\to Y$. Puis on prolonge $f$ au disque $D$ de sorte que la métrique induite sur $D$ soit à courbure négative ou nulle. Dans le cas où $Y$ est un polyèdre (le seul que nous traiterons ici), il suffit de prendre pour $f$ une surface réglée, par exemple un cône sur $f_{|\partial D}$, la métrique induite est alors riemannienne plate avec un nombre fini de singularités côniques. 

Suivant A. Weil \cite{W1}, on change la métrique induite $ds^{2}$ en une métrique riemannienne plate $ds_{0}^{2}\leq ds^{2}$ qui coïncide avec $ds^{2}$ le long du bord. En effet, par re\-pré\-sen\-ta\-tion conforme, quitte à composer $f$ avec un homéomorphisme quasiconforme du disque, $ds^{2}=e^{-2\phi}(du^{2}+dv^{2})$. La courbure $K=e^{2\phi}\Delta\phi$ étant négative ou nulle, la fonction harmonique $\phi_{0}$ qui vaut $\phi$ au bord satisfait $\phi\geq\phi_{0}$ à l'intérieur, et $ds_{0}^{2}=e^{-2\phi_{0}}(du^{2}+dv^{2})$ est plate. 

Suivant Yu. Reshetnyak, \cite{Reshetnyak1}, \cite{Reshetnyak2}, il existe un convexe $D'\subset\R^{2}$ et un ho\-mé\-o\-mor\-phis\-me $(D,ds_{0}^{2})\to D'$ qui augmente (au sens large) les distances, et préserve la longueur du bord. 

La restriction de cet homéomorphisme à $C_{k}$ est une application $g':C_{k}\to\R^{2}$ telle que, pour toute ar{\^e}te $cc'$ de $C_{k}$, $d(g'(c),g'(c'))=d(g(c),g(c'))$, et pour toute paire de sommets $c$, $c'$, $d(g'(c),g'(c'))\geq d(g(c),g(c'))$. Par cons{\'e}quent, $E_{1}(g')=E_{1}(g)$, et pour $j\geq 1$, $E_{j}(g')\geq E_{j}(g)$, donc
\begin{eqnarray*}
\frac{E_{1}(g)}{E_{j}(g)}\geq \frac{E_{1}(g')}{E_{j}(g')}\geq \frac{W(k,1)}{W(k,j)},
\end{eqnarray*}
d'apr{\`e}s la proposition \ref{rnwir}.\qed

\subsection{Bas du spectre, version Gromov}

Commençons par une variante de la définition \ref{defbas}. La terminologie est empruntée à \cite{Pichot}.

\begin{defi}
Soit $(C,m)$ un graphe pondéré. Soit $Y$ un espace métrique. Soit $g:C\to Y$. Le {\em quotient de Rayleigh, version Gromov} est $RQ^{Gro}(g)=E(g)/F(g)$, où
\begin{eqnarray*}
F(g)=\frac{1}{2\me}\sum_{c,\,c'}m(c)m(c')d(g(c),g(c'))^2 
\end{eqnarray*}
et $\me=\sum_{c}m(c)$.

Le {\em bas du spectre, version Gromov} $\lambda^{Gro}(C,m,Y)$ est la borne inférieure des $RQ^{Gro}(g)$, sur toutes les applications non constantes $g:C\to Y$.
\end{defi}

\begin{lemme}
\label{w>g}
Soit $Y$ un espace géodésique $CAT(0)$. Pour tout graphe pondéré $(C,m)$, 
$$\lambda^{Gro}(C,Y)\leq\lambda(C,Y).$$
L'égalité a lieu quand $Y$ est euclidien.
\end{lemme}

\preuve
Soit $c\in C^{0}$. Le lemme \ref{lembar} appliqué à $y=g(c)$ donne
\begin{eqnarray*}
\sum_{c'}m(c')d(g(c),g(c'))^2 \geq \sum_{c'}m(c')d(g(c'),bar(g))^2 +\me d(g(c),bar(g))^2 .
\end{eqnarray*}  
En sommant sur $c$, il vient
\begin{eqnarray*}
2\me F(g)&\geq& \sum_{c}m(c)d(g,bar(g))^2 +\me\sum_c m(c)d(g(c),bar(g))^2 \\
&=& 2\me d(g,bar(g))^2 ,
\end{eqnarray*}
avec égalité quand $Y$ est euclidien. On conclut que $\lambda^{Gro}\leq\lambda$, avec égalité quand $Y$ est euclidien. \qed

\begin{rem}
L'hypothèse $\lambda^{Gro}(liens,cones)>1/2$ sous laquelle M. Gromov prouve le théorème \ref{fixe} dans \cite{Gromov-RWRG} est donc un peu plus forte que celle de ce théorème. Au lieu d'utiliser le flot du gradient de l'énergie, M. Gromov en construit une approximation à temps discret, qui diminue strictement l'énergie. Pour des exposés de l'argument de \cite{Gromov-RWRG}, voir \cite{Pichot} ou \cite{Naor-Silberman}.
\end{rem}

\begin{rem}
Avec les notations du paragraphe \ref{sspec}, lorsqu'on choisit la pondération uniforme, pour un graphe dont chaque sommet a une valence au plus égale à $v$, $\displaystyle F\leq \frac{v^2}{\me}\sum_{j\geq 1}E_j$, et $\me$ est égal à deux fois le nombre d'arêtes de $C$. 
\end{rem}

\begin{lemme}
\label{wir=>lambda}
Si $Y$ satisfait à l'inégalité de Wirtinger $Wir_{k}$, alors le bas du spectre du cycle $C_{k}$ à valeurs dans $Y$ est au moins égal au bas du spectre à valeurs dans $\R$,
\begin{eqnarray*}
\lambda^{Gro}(C_{k},Y)\geq \lambda^{Gro}(C_{k},\R)=\frac{1}{2}|1-e^{2i\pi/k}|^{2}.
\end{eqnarray*}
\end{lemme}

\preuve
Si $g:C_{k}\to Y$,
\begin{eqnarray*}
RQ^{Gro}(g)^{-1}&=&
\frac{F(g)}{E(g)}\\
&=&\frac{1}{E(g)}\frac{1}{2k}\sum_{j\geq 1}E_{j}(g)\\
&\leq& \frac{1}{E(g)}\frac{1}{2k}\sum_{j\geq 1}\frac{W(k,j)}{W(k,1)}E_{1}(g)\\
&=&RG^{Gro}(g')^{-1},
\end{eqnarray*}
où $g'$ est un polygone régulier dans le plan euclidien.\qed

\subsection{\emph{Integralgeometrie}}

La formule de Garland, dans le cas des immeubles de type $\tilde{A}_2$, s'obtient en moyennant sur les appartements une formule naïve. Cette observation, développée dans \cite{Pansu2}, suggère d'estimer le bas du spectre des liens de ces immeubles en moyennant sur leurs appartements, qui sont des cycles.

\begin{prop}
\label{moyenne}
Soit $C$ un graphe de diamètre $\floor{k/2}$, dont tous les sommets ont une valeence $\leq v$, muni de la pondération uniforme. Soit $\mathcal{C}$ une famille non vide de $k$-cycles plongés isométriquement dans $C$. On suppose que, étant donnés deux sommets $c$ et $c'$ de $C$ le nombre $N_{\mathcal{C}}(c,c')$ de cycles de la famille $\mathcal{C}$ qui passent par $c$ et $c'$ ne dépend que de la distance $j=d(c,c')$, $N_{\mathcal{C}}(c,c')=N_{d(c,c')}=N_{j}$.

Soit $Y$ un espace métrique qui satisfait à l'inégalité de Wirtinger $Wir_{k}$. Alors le bas du spectre (au sens de Gromov) de $C$ à valeurs dans $Y$ est minoré par une constante 
\begin{eqnarray*}
\lambda^{Gro}(C,Y)\geq \lambda(k,v,N_j ),
\end{eqnarray*}
optimale au sens suivant. Si $C$ est régulier (tous les sommets ont la même valence) et admet un plongement $\iota$ dans un espace $Y_{0}$ tel que chaque cycle de $\mathcal{C}$ soit envoyé sur un polygone régulier contenu dans un plan euclidien plongé isométriquement dans $Y_{0}$, alors $\lambda(k)$ est égale au quotient de Rayleigh (selon Gromov) de ce plongement,
\begin{eqnarray*}
\lambda^{Gro}(C,Y)=RQ^{Gro}(\iota).
\end{eqnarray*}
\end{prop}

\preuve
Soit $\mathcal{PC}$ l'ensemble des paramétrisations de cycles de $\mathcal{C}$, i.e. des i\-so\-mé\-tries $h:C_{k}\to C$ dont l'image appartient à $\mathcal{C}$. Si $c$ et $c'\in C$ sont à distance $j\geq 1$, le nombre de triplets $(h,z,z')\in\mathcal{PC}\times C_{k}\times C_{k}$ tels que $h(z)=c$ et $h(z')=c'$ est égal à $N_{j}$ si $j\not=k/2$, et à $2N(k/2)$ sinon. En effet, pour tout couple $z$ et $z'$ à distance $j<k/2$, et pour tout cycle de $\mathcal{C}$ passant par $c$ et $c'$, il existe une unique paramétrisation qui envoie $z$ sur $c$ et $z'$ sur $c'$, et il en existe deux si $d(z,z')=k/2$. On notera donc
\begin{eqnarray*}
{\tilde N}_{j}=\left\{\begin{array}{ccc}
N_{j} \textrm{ si }j<k/2,\\
2N(k/2) \textrm{ sinon}.
\end{array}\right.
\end{eqnarray*}

Soit $g:C\to Y$ et $j\in\{1,\ldots,\floor{k/2}\}$,
\begin{eqnarray*}
\sum_{h\in\mathcal{PC}}E_{j}(g\circ h)
&=&\frac{1}{2}\sum_{h\in\mathcal{PC}}\sum_{z\in C_{k}}\sum_{\{z'\in C_{k}\,|\,d(z,z')=j\}} d(g\circ h(z),g\circ h(z'))^{2}\\
&=&\frac{1}{2}\sum_{c\in C}\sum_{\{c'\in C\,|\,d(c,c')=j\}} \sum_{\{(h,z,z')\,|\,h(z)=c,\,h(z')=c'\}}d(g(c),g(c'))^{2}\\
&=&
{\tilde N}_{j}E_{j}(g).
\end{eqnarray*}
Comme $Y$ satisfait $Wir_k$, pour tout $h$, $\displaystyle E_{j}(g\circ h)\leq \frac{W(k,j)}{W(k,1)}E_{1}(g\circ h)$, donc
\begin{eqnarray*}
E_{j}\leq \frac{W(k,j)}{{\tilde N}_{j}W(k,1)}E_{1}(g),
\end{eqnarray*}
\begin{eqnarray*}
F(g)\leq \frac{v^2}{\me}\sum_{j=1}^{\floor{k/2}}E_j (g)\leq E(g)\frac{v^2}{\me}\sum_{j=1}^{\floor{k/2}}\frac{W(k,j)}{{\tilde N}_{j}W(k,1)},
\end{eqnarray*}
soit $RQ^{Gro}(g)\geq \lambda$, où $\lambda$ est l'inverse de $\frac{v^2}{\me}\sum_{j=1}^{\floor{k/2}}\frac{W(k,j)}{{\tilde N}_{j}W(k,1)}$. 

Supposons qu'il existe un espace $CAT(0)$ $Y_0$ et un point $y_0 \in Y_0$ tels que, pour chaque cycle $\kappa$ de $\mathcal{C}$, $g_{|\kappa}$ factorise par le plongement de $\kappa$ comme polygone régulier d'un plan euclidien, puis par un plongement isométrique de ce plan dans $Y_0$ qui envoie le centre du polygone sur $y_0$. Alors l'inégalité $Wir_k$ est une égalité pour chaque $k$-cycle, donc $RQ^{Gro}(g)=\lambda(k)$.\qed

\begin{exemple}
Cas des immeubles euclidiens. 
\end{exemple}
Soit $I$ un immeuble euclidien localement fini de dimension 2. Soit $L$ le lien d'un sommet dans $I$. Alors $L$ est un immeuble sphérique de rang 1. Soit $\mathcal{C}$ l'ensemble de ses appartements. Chaque appartement est un $k$-cycle, avec $k=4$, $6$, $8$ ou $12$, il est plongé isométriquement. Le nombre d'appartements contenant deux points ne dépend que de leur distance. La proposition \ref{moyenne} s'applique, le cas d'égalité est réalisé par le plongement naturel du lien $\iota:L\subset I$.

\begin{cor}
\label{imm=>lambda>=}
Soit $I_p$ l'immeuble de Bruhat-Tits associé à $Sl(3,\Q_p )$, soit $L$ le lien d'un sommet de $I_p$, muni de la pondération uniforme. Soit $Y$ un espace métrique qui satisfait à l'inégalité $Wir_{6}$. Alors 
$$
\lambda^{Gro}(L,Y)\geq \lambda^{Gro}(L,I_p )=RQ^{Gro}(\iota:L\to I)=\frac{2(p^2 +p+1)}{7p^2 +4p+1}.
$$ 
\end{cor}

\preuve
$L$ est le graphe d'incidence du plan projectif sur le corps premier $\mathbb{F}_p$. Tous les sommets ont même degré $d=p+1$ (nombre de droites passant par un point donné, ou nombre de points sur une droite donnée), donc $m(c)=d$ pour tout $c\in L^0$. Notons $n$ le nombre de sommets.

Pour le plongement $\iota:L\subset I$, deux sommets voisins dans $L$ se trouvent à distance $1$ dans $I_p$, donc $E(\iota)=\frac{1}{2}nd$. Le nombre de couples de voisins est $nd$.

Comme $L$ n'a pas de cycles de longueur $<6$, le nombre de couples de sommets à distance 2 dans $L$ est égal au nombre de segments de longueur 2, soit $nd(d-1)$. Leur distance dans $I_p$ vaut $\sqrt{3}$.

Comme le diamètre de $L$ vaut 3, par deux points à distance 3 passent autant de chemins minimisants qu'un sommet possède de voisins, soit $d$. Le nombre de segments de longueur 3 vaut $nd(d-1)^2$, donc le nombre de couples de points à distance 3 vaut $nd(d-1)^2 /d=n(d-1)^2$. Leur distance dans $I_p$ vaut 2.

Il vient
\begin{eqnarray*}
F(\iota)&=&\frac{1}{2\me}\sum_{(c,c')}m(c)m(c')d(\iota(c),\iota(c'))^2 \\
&=&\frac{1}{2nd}d^2 (nd+3nd(d-1)+4n(d-1)^2 )\\
&=&\frac{d}{2}(1+(d-1)+3d(d-1)+4(d-1)^2 )\\
&=&\frac{1+p}{2}(1+4p+7p^2 ),
\end{eqnarray*}
d'où
\begin{eqnarray*}
RQ^{Gro}(\iota)=\frac{E(\iota)}{F(\iota)}=\frac{2(1+q+q^2 )}{1+4p+7p^2 }.\qed
\end{eqnarray*}

Malheureusement, $\lambda^{Gro}(L,I_p )<\frac{1}{2}$ pour tout $p\geq 2$, donc aucun théorème de point fixe n'en résulte. Noter que $RQ(\iota)=\frac{1}{2}$, et il est vraisemblable que $\lambda(L,Y)=\frac{1}{2}$ pour tout espace géodésique $CAT(0)$ $Y$.

\subsection{Spéculation}

Voici une version simplifiée de la Proposition \ref{moyenne}.

\begin{prop}
Soit $Y$ un espace géodésique $CAT(0)$. Soit $C$ un graphe fini. Soit $\mathcal{L}$ une famille de lacets de longueur $\leq k$ (i.e. des applications du $j$-cycle, $j\leq k$, dans $C$). On note
\begin{itemize}
  \item $v$ la borne supérieure des valences des sommets de $C$.
  \item $r$ la borne inférieure du nombre de lacets dans $\mathcal{L}$ qui contiennent une paire de sommets de $C$ donnée.
  \item $q$ la borne supérieure du nombre de lacets dans $\mathcal{L}$ qui contiennent une arête donnée.
  \item $A$ le nombre d'arêtes de $C$.
\end{itemize}
Alors
\begin{eqnarray*}
\lambda(C,Y)\geq \lambda^{Gro}(C,Y)\geq\frac{4Ar}{qkv^2}\frac{1}{2}|1-e^{2i\pi/k}|^2 .
\end{eqnarray*}
En particulier, si $k=6$ et $\frac{Ar}{qv^2}>\frac{3}{2}$, alors $\lambda^{Gro}(C,Y)>\frac{1}{2}$.
\end{prop}

\preuve
On rappelle qu'on utilise la pondération uniforme, pour laquelle $m=1$ pour les arêtes, et pour un sommet, $m(c)$ est la valence de $c$. Par exemple, $m(z)=2$ pour tout sommet $z$ du $j$-cycle $C_j$. Chaque élément de $\mathcal{L}$ a une longueur $j(\ell)$, i.e. est une application du $j(\ell)$-cycle $C_{j(\ell)}$ dans $C^0$. Par définition de $r$,
\begin{eqnarray*}
N(c,c')=\textrm{card}(\{(\ell,z,z')\in\mathcal{L}\times C_{j(\ell)}\times C_{j(\ell)}\,|\,\ell(z)=c,\,\ell(z')=c'\})\geq r.
\end{eqnarray*}
Il vient
\begin{eqnarray*}
\sum_{\ell\in\mathcal{L}}F(g\circ\ell)&=&\sum_{\ell}\frac{1}{4j(\ell)}\sum_{z,\,z'\in C_{j(\ell)}} m(z)m(z')d(g\circ\ell(z),g\circ\ell(z'))^2 \\
&\geq&\frac{1}{k}\sum_{\ell}\sum_{z,\,z'\in C_{j(\ell)}}d(g\circ\ell(z),g\circ\ell(z'))^2 \\
&\geq&\frac{1}{kv^2}\sum_{c,\,c'\in C}v^2 N(c,c')d(g(c),g(c'))^2 \\
&\geq&\frac{r}{kv^2}\sum_{c,\,c'\in C} m(c)m(c')d(g(c),g(c'))^2 =\frac{4Ar}{kv^2}F(g).
\end{eqnarray*}
On applique le Lemme \ref{wir=>lambda} à chaque lacet $\ell\in\mathcal{L}$,
\begin{eqnarray*}
E(g\circ\ell)&\geq& \lambda^{Gro}(C_{j(\ell)},Y)F(g\circ\ell)\\
&\geq&\frac{1}{2}|1-e^{2i\pi/k}|^{2}F(g\circ\ell).
\end{eqnarray*}
Il vient
\begin{eqnarray*}
\frac{1}{2}|1-e^{2i\pi/k}|^{2}\sum_{\ell\in\mathcal{L}}F(g\circ\ell)
&\leq&\sum_{\ell\in\mathcal{L}}E(g\circ\ell)\\
&=&\sum_{\ell\in\mathcal{L}}\sum_{e\textrm{ edge of }C_{j(\ell)}}d(g\circ\ell(ori(e)),g\circ\ell(extr(e)))^2 \\
&\leq&q\sum_{e'\textrm{ edge of }C}d(g(orig(e')),g(extr(e'))^2 \\
&=&q\,E(g),
\end{eqnarray*}
et enfin
\begin{eqnarray*}
RQ(g)\geq RQ^{Gro}(g)\geq \frac{4Ar}{qkv^2}\frac{1}{2}|1-e^{2i\pi/k}|^2 .\qed
\end{eqnarray*}

\medskip

\textbf{Question}. 
Un graphe aléatoire à $n$ sommets, de degré $n^{1/3}<d<n^{1/2}$ admet-il une famille de lacets de longueur $\leq 6$ qui satisfait $\frac{Ar}{qv^2}>\frac{3}{2}$ ?

\section{Modèles de groupes aléatoires à densité}

\subsection{Les modèles}

D'après Gromov \cite{Gromov-AIIG}, un modèle de groupes aléatoires consiste à se donner une famille de lois de probabilité sur l'ensemble des présentations de groupes. Une présentation de groupe, c'est un ensemble de générateurs $S$ et un ensemble de relateurs $R$, i.e. un sous-ensemble du groupe libre $F_{S}$. Dans ce texte, on va décrire des modèles portant sur des présentations finies. Pour simplifier, on se limite à des lois uniformes sur des ensembles finis. Il y a alors trois paramètres en jeu : le nombre $m$ d'éléments de $S$, la longueur (ou bien longueur maximale) $\ell$ des éléments de $R$ et le nombre $N$ de relateurs. Un trop grand nombre de relateurs risque de rendre le groupe quotient $G=\langle S\,|\,R\rangle$ trivial. Gromov a découvert que même en prenant pour $N$ une puissance du nombre de choix possibles, le groupe obtenu est souvent non trivial. Cela l'a conduit à choisir $N$ de la forme
\begin{eqnarray*}
N=\frac{2m}{2m-1}(2m-1)^{d\ell},
\end{eqnarray*}
où $0<d<1$. On parle alors de modèle {\em à densité} $d$.

Si on fixe $m$ et on fait tendre $\ell$ vers l'infini, on obtient le modèle étudié par Gromov dans \cite{Gromov-AIIG}. Si on fixe $\ell=3$ et on fait tendre $m$ vers l'infini, on obtient un modèle étudié par A. Zuk dans \cite{Zuk-03}. On s'intéresse à des énoncés du type suivant. \emph{On se place en densité $d$, i.e. on fixe $N\sim (2m-1)^{d\ell}$. La proportion des présentations possédant la propriété (P) tend vers 1 lorsque $\ell$ (resp. $m$) tend vers l'infini}. Lorsqu'une telle assertion est vraie, on dit que \emph{la propriété (P) est génériquement vraie en densité $d$}.

\begin{theo}
\label{gralinfini}
\emph{(Gromov, \cite{Gromov-AIIG}, \cite{Ollivier}).} En densité $d>1/2$, le groupe $G=\langle S\,|\,R\rangle$ est génériquement trivial. En densité $d<1/2$, $G$ est génériquement hyperbolique et infini.
\end{theo}

\subsection{Généricité de la propriété (T)}

A une présentation finie d'un groupe est associé un polyèdre fini $\mathcal{P}_{S,R}$, le \emph{polyèdre de Cayley}, défini comme suit. On forme un bouquet de cercles, un pour chaque générateur $s\in S$. Pour chaque relateur $r\in R$, on colle une 2-cellule sur le bouquet au moyen de l'application qui paramètre par le bord d'un disque le chemin décrit par le mot $r$. Par construction, le groupe fondamental du polyèdre de Cayley d'une présentation de $G$ est $G$. Par conséquent, son revêtement universel comporte autant de sommets qu'il y a d'éléments de $G$. On peut construire directement le revêtement universel $\tilde{\mathcal{P}}_{S,R}$ comme suit. Il y a une arête numérotée $s$ reliant deux éléments $g$ et $g'$ de $G$ si $g'=gs$. Pour chaque relateur $r$ et chaque élément $g\in G$, on colle une 2-cellule le long du chemin issu de $g$ obtenu en suivant successivement les arêtes dont les numéros sont indiqués par le mot $r$.

Lorsque les relateurs sont de longueur 3, les 2-cellules sont des triangles. Toutefois, $\mathcal{P}_{S,R}$ n'est pas un complexe simplicial, puisqu'il n'a qu'un sommet. Son re\-vê\-te\-ment universel $\tilde{\mathcal{P}}_{S,R}$ non plus, car il arrive que deux arêtes aient les mêmes extrémités, ou que deux faces soient collées sur le même chemin. Néanmoins, on peut définir le lien d'un sommet dans $\tilde{\mathcal{P}}_{S,R}$. C'est un graphe $L(S,R)$ dont l'ensemble des sommets s'identifie à $S\cup S^{-1}$. Il y a une arête reliant $z$ à $z'$ chaque fois qu'il existe $z"\in S\cup S^{-1}$ tel que $z^{-1}z'z''$ ou $z'^{-1}z''z$ ou $z''^{-1}zz'$ appartient à $R$.

\begin{theo}
\label{gralT}
\emph{(Zuk, \cite{Zuk-03}).} On se place dans le modèle où les relateurs sont de longueur 3. Si $d<1/3$, génériquement $G$ possède un quotient libre à deux gé\-né\-ra\-teurs. Si $d>1/3$, génériquement le graphe $L(S,R)$ a un bas du spectre proche de 1. En particulier, il satisfait
\begin{eqnarray*}
\lambda(L(S,R),\R)>\frac{1}{2}.
\end{eqnarray*}
\end{theo}

\preuve
On introduit un modèle de graphes aléatoires, et on montre que gé\-né\-ri\-que\-ment en densité $>1/3$, $L(S,R)$ a les mêmes propriétés spectrales qu'un graphe aléatoire.\qed

\begin{cor}
\label{corzuk}
Dans le modèle où les relateurs sont de longueur 3, et en densité $1/3<d<1/2$, le groupe associé à une présentation générique est infini et possède la propriété (T) de Kazhdan. En revanche, en densité $<1/3$, il ne possède pas la propriété (T).
\end{cor}

\preuve
Comme le polyèdre de Cayley n'est pas à proprement parler un complexe simplicial, le résultat ne découle pas littéralement des méthodes du paragraphe \ref{sfixe}, mais l'esprit est le même. Voir \cite{Zuk-03} pour les détails.\qed

\subsection{Autres modèles}

Dans \cite{Gromov-RWRG}, Gromov introduit un procédé de construction de groupes dont le graphe de Cayley est modelé sur un graphe donné. Il y a un choix supplémentaire, celui d'un étiquetage du groupe. Comme on peut faire ce choix au hasard, on obtient ainsi un nouveau modèle de groupe aléatoire. Pour des exposés de cette construction, voir \cite{Ghys} et \cite{Ollivier2}.

Des idées voisines de celles exposées dans ce texte (présence dans le graphe de Cayley d'un groupe d'un sous-graphe de grand $\lambda_1$) permettent de prouver que la propriété (T) (section 3 de \cite{Gromov-RWRG}, \cite{Silberman}) et une propriété de point fixe sur une classe d'espaces CAT(0) (\cite{Naor-Silberman}) sont génériques aussi dans ce modèle.

\subsection{Conclusion}

Je ne crois pas qu'on connaisse des groupes infinis qui possèdent la propriété de point fixe sur tous les espaces CAT(0). Il est tentant d'appliquer le théorème de superrigidité de N. Monod, \cite{Monod}, mais les conditions de non évanescence ou de réductivité en limitent pour l'instant le champ d'application. 

A ma connaissance, la question de savoir si la propriété de point fixe sur tous les espaces CAT(0) est générique dans un modèle de groupes aléatoires est encore ouverte. 

La propriété FSI est strictement plus faible. En effet, M. Burger et S. Mozes ont construit dans \cite{Burger-Mozes1}, \cite{Burger-Mozes2} des groupes simples (ils ont donc la propriété FSI) qui sont des réseaux dans des produits d'arbres, donc agissent sans point fixe sur un immeuble (voir aussi \cite{Caprace-Remy}). Pourtant, on ne sait pas plus si elle est générique.

\par\medskip\noindent
Pierre Pansu\\
Laboratoire de Mathématique d'Orsay\\
UMR 8628 du C.N.R.S.\\
Bâtiment 425\\
Université Paris Sud - 91405 Orsay (France)\\
\smallskip\noindent
{\tt\small Pierre.Pansu@math.u-psud.fr}\\
http://www.math.u-psud.fr/$\sim$pansu

\end{document}